\documentclass{m2an}
\usepackage{hyperref}
\usepackage{graphicx}
\usepackage{ulem}
\usepackage{siunitx}  
\allowdisplaybreaks

\usepackage[colorinlistoftodos]{todonotes}

\newcommand{\F}{\mathrm{f}}
\newcommand{\BJ}{\mathrm{BJ}}
\newcommand{\PM}{\mathrm{pm}}
\newcommand{\HYB}{\mathrm{hyb}}
\newcommand{\PS}{\mathrm{ps}}
\newcommand{\ten}[1]{\ensuremath{\boldsymbol{\mathsf{#1}}}}
\renewcommand{\vec}[1]{{\ensuremath{\boldsymbol{\mathrm #1}}}}
\renewcommand{\d}{\mathrm{d}}
\newcommand{\Z}{L^2(\Gamma)}
\newcommand{\Ce}{C_{\mathrm{emb}}}
\newcommand{\MI}{M_{1,1}^{1,\mathrm{bl}}}
\newcommand{\MII}{M_{2,1}^{1,\mathrm{bl}}}

\newtheorem{theorem}{Theorem}
\newtheorem{lemma}{Lemma}

\newtheorem{remark}{Remark}

\definecolor{MidnightBlue}{cmyk}{0.98,0.13,0,0.43}

\begin{document}
%%-----------------------------
%%      the top matter
%%-----------------------------
\title{An advanced hybrid-dimensional Stokes--Darcy model for fractured porous media}
\thanks{The first and the third author thank the Deutsche Forschungsgemeinschaft (DFG, German Research Foundation) -- Project Number 327154368 -- SFB 1313 for the financial support. The second and the third author thank the DFG -- Project Number 490872182 for the financial support.}

\author{Paula Strohbeck}\address{Institute of Applied Analysis and Numerical Simulation, University of Stuttgart, Stuttgart, Germany.}
\author{Linheng Ruan}\sameaddress{1}
\author{Iryna Rybak}\sameaddress{1}
\date{\today}
\begin{abstract}
%Title should be simple and informative. A shortened version of the title consisting of a maximum of 75 characters (including spaces) for running headers should also be provided. 
%An abstract in English is required. It should be completely self-contained, not exceeding 200 words and written as a single paragraph.

Thin fluid-filled channels in porous media are prevalent in a wide range of environmental settings including fractures in geological formations. Transport of fluid in such multi-domain systems is commonly described by Darcy’s law 
and formulated either in a fully dimensional framework or through a hybrid-dimensional approach, where fractures are represented as inclusions of co-dimension one enclosed in the porous matrix. These mixed-dimensional models enable accurate  representation of the corresponding full-dimensional systems while significantly reducing computational complexity. Models available in literature usually do not properly account for channel-wall roughness while applying coupling conditions on the fracture-matrix interface and are usually developed for flow directions parallel to the fracture. The goal of this work is to advance and study the hybrid-dimensional Stokes--Darcy model developed by the authors in their previous work. We consider a thin open channel embedded within a porous matrix, where the channel walls are in contact with porous regions of distinct  permeabilities. We incorporate channel-wall roughness effects into the model and improve it to accommodate arbitrary flow directions. We prove existence and uniqueness of the weak solution for the new model and demonstrate its advantages over existing hybrid-dimensional Stokes--Darcy models using pore-scale resolved simulations.

\end{abstract}

\maketitle
%%-----------------------------

\section*{Introduction}

Thin fluid-filled channels embedded in porous media appear in a wide range of environmental, industrial and medical applications including fractures in geological formations, conduits in karst aquifers, vascular networks in biological tissues, and engineered flow channels in reactors and heat exchangers. 
Models for fractured porous media can generally be classified into two categories: continuous and discrete fracture-matrix models.

In the first approach, fractures are represented implicitly, and the fractured porous medium is treated as a continuum in which the effects of fractures are incorporated into effective properties, such as permeability and porosity, e.g.,~\cite{Arbogast1990}. In the second approach, fractures are represented explicitly as inclusions embedded in the surrounding porous matrix, e.g.~\cite{Martin_etal_05, Ahmed_etal_17, frih2008modeling, Flemisch_etal_18, AngotHubert2009, Berre2019, Mikelic2015, rivas2025fluid}. 
The continuum approach is computationally attractive because it avoids explicitly resolving individual fractures. However, its applicability relies on the assumption that the fracture network can be adequately homogenised. This assumption may not be valid when fractures are highly connected or characterised by the length scale comparable to those of the domain of interest. Under such conditions, conventional homogenisation  fails to capture the preferential flow paths. In contrast, discrete fracture models explicitly represent individual fractures embedded in the porous matrix. Such models preserve important geometrical characteristics of fracture networks, including fracture location, orientation, aperture, length, and connectivity. This explicit representation is particularly important in fractured porous media, where fluid transport may be dominated by preferential pathways formed by interconnected fractures. Discrete formulations also provide a framework for describing flow within fractures and fluid exchange between the fracture and the surrounding porous matrix. In this work, we are dealing with the discrete fracture modelling.  

Transport of fluid through such  multi-domain systems is typically described by Darcy's law or its extensions, both in the fracture and in the surrounding porous matrix, applying appropriate coupling conditions at the boundary between two subdomains
(fracture-matrix interface). Such flow models can be formulated either in a full-dimensional setting, where the fractures and the surrounding porous media are represented in the same spatial dimension, or using a hybrid-dimensional framework, in which fractures are modelled as lower-dimensional interfaces of co-dimension one embedded in the porous matrix. When properly derived, such hybrid-dimensional models provide accurate approximations of the corresponding full-dimensional formulations while significantly reducing computational costs, as they eliminate the need for highly refined meshes within the fracture domain. 

There exist several mathematical techniques to derive hybrid-dimensional models for fractured porous media: vertical averaging, e.g.~\cite{Knabner_Roberts_14, Formaggia_etal_14, Lesinigo_etal_11, Rybak-Metzger-20, Martin_etal_05,  brenner-etal-2018, Bukac-et-al-17, BoonNordbotten2023, Berre2019, BudisaHu2021, Ahmed_etal_17, frih2008modeling, Flemisch_etal_18, Glaeser_etal_17}, asymptotic modelling~\cite{AngotHubert2009}, homogenisation theory~\cite{maxi2024, List2020, morales2017darcy,Morales2012} and dimensional reduction by Fourier analysis~\cite{gander-hennicker-masson-21}. For these multi-domain flow systems,
different mathematical models can be used in the fracture and in the porous matrix: Darcy--Darcy, e.g.~\cite{Martin_etal_05, AngotHubert2009,Formaggia_etal_14, Burbulla2023,maxi2024, gander2023dimensional}, Brinkman--Darcy~\cite{Lesinigo_etal_11,Morales2017}, Darcy--Forchheimer~\cite{frih2008modeling, Knabner_Roberts_14} and Stokes--Darcy~\cite{Rybak-Metzger-20,gander-hennicker-masson-21} coupled models. For fluid-filled channels without any debris, Stokes--Darcy models are the most suitable choice. For two-fluid-phase flows in fractured porous media, models are based on two-phase Darcy's law, e.g.~\cite{Ahmed_etal_17, brenner-etal-2018, BoonNordbotten2023}. Here, we provided a brief overview of mathematical models for fractured porous media. For a comprehensive review on modelling concepts and efficient solution strategies, we refer the reader to~\cite{Berre2019}.

 In this work, we deal with single-phase flows and enhance the Stokes--Darcy model presented in our previous work~\cite{Rybak-Metzger-20}. We derive the advanced hybrid-dimensional model by means of vertical averaging. The Stokes equations are used in the fracture domain and averaged across the fracture in the normal direction. To get a closed model formulation, a suitable set of coupling conditions at the fracture-matrix interface is needed. In contrast to~\cite{Rybak-Metzger-20}, where the classical set of conditions containing the Beavers--Joseph or Beavers--Joseph--Saffman condition on the tangential velocity component was considered, we now apply the generalised interface conditions developed in~\cite{Eggenweiler_Rybak_MMS20}.
They generalise the classical coupling conditions by incorporating the pore-scale surface roughness of the fracture-matrix interface through rigorously determined coefficients computed using homogenisation and boundary layer theory. These generalised conditions are applicable to arbitrary flow directions~\cite{Eggenweiler_Rybak_20}.

The goals of this work are (i) development of an advanced hybrid-dimensional Stokes--Darcy model using the generalised coupling conditions, (ii) well-posedness analysis for the new model and (iii) its validation against the pore-scale resolved model and other hybrid-dimensional Stokes--Darcy models available in the literature. The paper is structured as follows. In Section~\ref{sec:geometry-assumptions}, we present the geometrical setting and the assumptions used to develop the hybrid-dimensional model. In Section~\ref{sec:fracture-models}, we provide two classes of fracture models, including the pore-scale resolved model and the macroscale full-dimensional models, which serve as the basis for development of hybrid-dimensional models. Section~\ref{sec:Reduced-dim-model} is devoted to the derivation of the advanced hybrid-dimensional Stokes--Darcy model. The well-posedness of this proposed model is proved in Section~\ref{sec:ana}. Due to additional terms arising in the generalised interface conditions, the proof is more challenging than in the case of Beavers--Joseph--Saffman condition~\cite{Rybak-Metzger-20}. In Section~\ref{sec:NumSimulation}, we validate the newly developed model against the pore-scale resolved model and compare it to the previously developed models~\cite{Rybak-Metzger-20} to demonstrate its advantages. Conclusions and future work follow in Section~\ref{sec:Conclusion}.

\section{Modelling assumptions and geometrical setting}
\label{sec:geometry-assumptions}
\begin{figure}[!ht]
    \centering
    \raisebox{0.5mm}{
    \includegraphics[height=2.8cm]{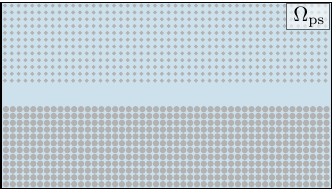}} 
    \hfill
    \includegraphics[height=2.8cm]{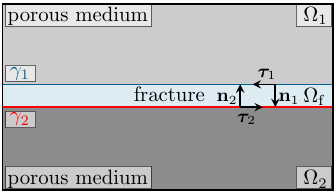} 
    \hfill %\hspace*{5ex}
    \includegraphics[height=2.8cm]{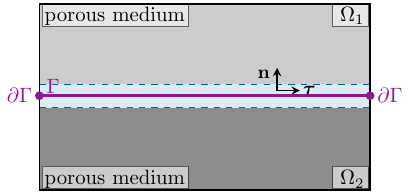}
    \caption{Geometry of the fracture models: pore-scale (left), full-dimensional macroscale (middle) and hybrid-dimensional macroscale (right).}
    \label{fig:fracture-model}
\end{figure}

In this work, we consider three different models for an open channel embedded in two distinct porous media: (i) the pore-scale resolved model (Fig.~\ref{fig:fracture-model}, left), (ii) the macroscale full-dimensional model (Fig.~\ref{fig:fracture-model}, middle), and (iii)~the macroscale hybrid-dimensional model (Fig.~\ref{fig:fracture-model}, right). The full-dimensional model will be applied to derive the dimensionally reduced model for the fracture and, consequently, for the overall hybrid-dimensional model. 
The pore-scale resolved model will be used for validation purposes to compare the new hybrid-dimensional model with two other presented in~\cite{Rybak-Metzger-20}.

In this work, we restrict ourselves to a Lipschitz domain in two space dimensions $\overline{\Omega} = \overline{\Omega}_1 \cup \overline{\Omega}_\F \cup \overline{\Omega}_2\subset\mathbb{R}^2$ composed of two porous-medium regions $\Omega_1,\, \Omega_2\subset\mathbb{R}^2$ divided by the fracture $\Omega_\F\subset\mathbb{R}^2$ (Fig.~\ref{fig:fracture-model}). The fracture $\Omega_\F$ is assumed to be horizontal and of constant thickness $d > 0$. 
The sharp interfaces on the top of the fracture, $\gamma_1=\overline{\Omega}_1 \cap \overline{\Omega}_\F \setminus \partial \Omega$, and on the bottom, $\gamma_2=\overline{\Omega}_2\cap \overline{\Omega}_\F \setminus \partial \Omega$, are considered to be straight. The pore-scale surface roughness at the fracture-matrix interfaces is already incorporated in the generalised interface conditions. 

At the pore-scale, the porous media are composed of periodically distributed solid inclusions, which allows determination of the model parameters such as permeability values and boundary layer constants. We assume separation of scales $\varepsilon = \ell / L \ll 1$, where $\varepsilon$ is the scale separation parameter, $\ell$ is the characteristic pore size and $L$ is the length of the flow domain. Note that this geometrical assumption is not a restriction on the derivation of the hybrid-dimensional model, which is also applicable to heterogeneous and non-periodic porous media.

We assume a single-phase steady-state flow at low Reynolds numbers $Re \ll 1$. The fluid is incompressible and the viscosity $\mu > 0$ is constant. The flow system is assumed to be isothermal. The porous media are rigid, non-deformable and fully saturated with the same fluid as in the fracture.

 \section{Fracture models}\label{sec:Full-dim-model}
 \label{sec:fracture-models}
In this section, we introduce the mathematical models for flow in fractured porous media, both at pore-scale and macroscale.

\subsection{Pore-scale resolved model}
The flow in the fluid domain $\Omega_{\PS}$ is described by the Stokes equations with the no-slip condition on the boundaries of the solid inclusions
\begin{subequations}
\label{eq:Stokes-micro}
\begin{align}
    \nabla \cdot \vec{v}_{\PS} = 0, \quad - \nabla \cdot \ten{T}(\vec{v}_{\PS}, p_{\PS}) = \vec{f}_{\PS} \qquad &\text{in} \; \; \Omega_{\PS},\\
     \vec{v}_{\PS} = \vec{0} \qquad &\text{on} \; \; \partial \Omega_{\PS}\backslash \partial \Omega,
\end{align}
\end{subequations}
with the pore-scale flow velocity $\vec{v}_{\PS}$, the pore-scale pressure $p_{\PS}$, the source term $\vec{f}_{\PS}$ and 
$\ten{T}(\vec{v}_{\PS}, p_{\PS}) := \mu \nabla \vec{v}_{\PS} - p_{\PS} \ten{I}$ the stress tensor with dynamic viscosity $\mu > 0$. On the external boundary $\partial \Omega = \Gamma_D \cup \Gamma_N, \; \Gamma_D \cap \Gamma_N = \emptyset$, we impose Dirichlet or Neumann boundary conditions
\begin{align}
\label{eq:BC-micro}
    \vec{v}_{\PS} = \overline{\vec{v}}_\PS \; \; \text{on} \;\; \Gamma_D, \qquad \ten{T}\left(\vec{v}_{\PS},p_{\PS}\right)\vec{n}= \overline{\vec{h}}_\PS \;\; \text{on} \;\; \Gamma_N,
\end{align}
with given functions $\overline{\vec{v}}_\PS$ and $\overline{\vec{h}}_\PS$, and the unit outward normal vector $\vec{n}$ on $\partial \Omega$.

\subsection{Full-dimensional Stokes--Darcy models}
In the full-dimensional macroscale models, the fluid flow in the porous-medium domains $\Omega_i,\; i \in \{1,2\}$ is described by Darcy's law
\begin{subequations}
\label{eq:Darcy-law}
\begin{align}
    \nabla \cdot \vec{v}_{i} = f_{i}, \quad \vec{v}_{i} = - &\frac{\ten{K}_i}{\mu}\nabla p_{i} &&\text{in} \; \; \Omega_{i}, \label{eq:Darcy-dual}\\
\intertext{where $\vec{v}_{i}$ is the Darcy velocity, $p_{i}$ is the pressure, $f_i$ is a source term and $\ten{K}_i \in \mathbb{R}^{2\times 2}$ is the permeability of the respective domain $\Omega_i$. In this work, we use the porous-medium model in its primal form}
    -\nabla \cdot \left(\frac{\ten{K}_i}{\mu}\nabla p_{i}\right) &= f_i &&\text{in} \; \; \Omega_{i}. \label{eq:Darcy-primal}
\end{align}
\end{subequations}
Both permeability tensors are symmetric, positive definite and bounded
\begin{eqnarray}
    k_{\min,i} \| \vec{x} \|^2  \le \vec{x} \cdot \ten{K}_i \vec{x} \le k_{\max,i} \| \vec{x} \|^2, \quad \forall \vec{x} \in \mathbb{R}^2, \quad i\in\{1,2\},\label{equ:boundedK}
\end{eqnarray}
where $  k_{\max,i}\ge k_{\min,i}>0$ for $i\in\{1,2\}$.
We set Dirichlet or Neumann boundary conditions on the external boundary of the porous-medium domains $\Gamma_i = \partial \Omega_i \backslash \gamma_i$:
\begin{align}
\label{eq:BC-Darcy}
    p_i = \overline{p}_i \text{ on } \Gamma_{D,i}, \quad \vec{v}_i \cdot \vec{n}_i = \overline{v}_i \text{ on } \Gamma_{N,i}, \quad i=1,2,
\end{align}
where $\overline{p}_i$ and $\overline{v}_i$ are given and $\Gamma_i = \Gamma_{D,i}\cup \Gamma_{N,i}, \; \Gamma_{D,i}\cap \Gamma_{N,i} = \emptyset, \; \Gamma_{D,i} \neq \emptyset$. 

The flow in the fracture $\Omega_\F$ is governed by the Stokes equations
\begin{subequations}
\label{eq:Stokes-macro}
    \begin{align}
        \nabla \cdot \vec{v}_\F = 0, \quad &\text{in} \; \; \Omega_\F,\label{eq:Stokes1}\\
         - \nabla \cdot \ten{T}(\vec{v}_\F, p_\F) = \vec{f}_\F \quad &\text{in} \; \; \Omega_\F,\label{eq:Stokes2}
\end{align}
\end{subequations}
where $\vec{v}_\F$ is the fluid velocity, $p_\F$ is the fluid pressure, $\ten{T}(\vec{v}_\F, p_\F)$ is the stress tensor and $\vec{f}_\F$ is a source term. On the external boundary of the fracture $\Gamma_\F = \partial \Omega_\F \backslash \{\gamma_1 \cup \gamma_2\}$, we have the boundary conditions
\begin{equation}
\label{eq:BC-Stokes}
    \vec{v}_\F = \overline{\vec{v}}_\F \text{ on } \Gamma_{D,\F}, \quad \ten{T}(\vec{v}_\F, p_\F) \vec{n} = \overline{\vec{h}}_\F \text{ on } \Gamma_{N,\F},
\end{equation}
where $\overline{\vec{v}}_\F$ and $\overline{\vec{h}}_\F$ are given functions and $\Gamma_\F = \Gamma_{D,\F} \cup \Gamma_{N, \F}$ and $\Gamma_{D,\F} \not = \emptyset$. 
Note that $\Gamma_\F$ is composed of two different segments. As we average along the segments to get boundary conditions for the hybrid-dimensional model, mixed Dirichlet/Neumann conditions are not allowed on the same segment.

To complete the full-dimensional macroscale formulations, which provide the basis for the hybrid-dimensional models, the Stokes--Darcy systems~\eqref{eq:Darcy-law}--\eqref{eq:BC-Stokes} need to be completed by appropriate coupling conditions at the fracture-matrix interfaces $\gamma_1$ and $\gamma_2$ (Fig.~\ref{fig:fracture-model}, middle). 

In the literature, there exist several possibilities for the interface conditions between the Stokes and Darcy regions. Depending on the chosen coupling conditions, we obtain different hybrid-dimensional models. The classical set of coupling conditions consists of the conservation of mass~\eqref{eq:IC-MC}, the balance of normal forces~\eqref{eq:IC-BONF} and either the Beavers--Joseph~\eqref{eq:IC-BJ} or the Beavers--Joseph--Saffman~\eqref{eq:IC-BJS} condition on the tangential velocity
\begin{subequations}
\label{eq:IC-classical}
\begin{align}
    \vec{v}_\F \cdot \vec{n}_i &=\vec{v}_i \cdot \vec{n}_i \quad &&\text { on } \gamma_i, \label{eq:IC-MC}\\
    -\vec{n}_i \cdot \ten{T}\left(\vec{v}_\F, p_\F\right) \vec{n}_i   &= p_i \quad&&\text { on } \gamma_i, \label{eq:IC-BONF}\\
    \left(\vec v_\F - \vec v_i\right) \cdot \boldsymbol{\tau}_i - \alpha_i^{-1} \sqrt{K}_i \left(\nabla \vec v_\F \cdot \vec{n}_i\right) \boldsymbol{\tau}_i &= 0 \quad &&\text { on } \gamma_i, \label{eq:IC-BJ} \\
    \vec{v}_\F \cdot \boldsymbol{\tau}_i - \alpha_i^{-1} \sqrt{K}_i \left(\nabla \vec v_\F \cdot \vec{n}_i\right) \boldsymbol{\tau}_i &= 0 \quad &&\text { on } \gamma_i. \label{eq:IC-BJS}
\end{align}
\end{subequations}
Here, $\vec n_i$ is the unit normal vector pointing outward from the porous medium $\Omega_i$ and $\boldsymbol{\tau}_i$ is the corresponding tangential vector (Fig.~\ref{fig:fracture-model}, middle). The following relations hold for the normal and tangential vectors at the lower dimensional interface in the hybrid-dimensional model: $\vec{n} = -\vec{n}_1 = \vec{n}_2$ and $\boldsymbol{\tau} = -\boldsymbol{\tau}_1 = \boldsymbol{\tau}_2$ (Fig.~\ref{fig:fracture-model}, right). The constant $\alpha_i > 0$ is the Beavers-Joseph slip coefficient.

Several approaches have been proposed in the literature for computing $\sqrt{K_i}$. An overview can be found in~\cite[Table 1]{Eggenweiler_Rybak_20}. In this work, we choose $\sqrt{K_i} = \sqrt{\boldsymbol{\tau}_i \cdot \ten{K}_i \boldsymbol{\tau}_i}$, which is most commonly used. 

To derive the advanced hybrid-dimensional model, we apply the generalised coupling conditions from our previous work~\cite{Eggenweiler_Rybak_MMS20,Strohbeck-Eggenweiler-Rybak-23} instead of the classical conditions given in~\eqref{eq:IC-classical}. They comprise the conservation of mass~\eqref{eq:IC-ER-MC}, an extension of the balance of normal forces~\eqref{eq:IC-ER-BONF} and a generalisation of the Beavers--Joseph condition~\eqref{eq:IC-ER-BJ}:
\begin{subequations}
\label{eq:IC-ER}
\begin{align}
    \vec{v}_\F \cdot \vec{n}_i &=\vec{v}_i \cdot \vec{n}_i \quad &&\text { on } \gamma_i, \label{eq:IC-ER-MC}\\
    -\vec{n}_i \cdot \ten{T}\left(\vec{v}_\F, p_\F\right) \vec{n}_i - \mu N_{i,s}^{\mathrm{bl}} \; \boldsymbol{\tau}_i \cdot \nabla \vec{v}_\F \,\mathbf{n}_i  &= p_i \quad&&\text { on } \gamma_i, \label{eq:IC-ER-BONF}\\
    \vec v_\F \cdot \boldsymbol{\tau}_i - \ell N_{i,1}^{\mathrm{bl}} \,\boldsymbol{\tau}_i \cdot \ten T(\vec v_\F, p_\F) \vec n_i &= {-} \frac{\ell^2}{\mu} ( \ten{M}_i^{\mathrm{bl}} \nabla p_i ) \cdot \vec{\tau}_i  \quad &&\text { on } \gamma_i. \label{eq:IC-ER-BJ}
\end{align}
\end{subequations}
The model parameters $N_{i,s}^{\mathrm{bl}} \in \mathbb{R},\, N_{i,1}^{\mathrm{bl}} \in \mathbb{R}$ and $\ten{M}_i^{\mathrm{bl}} \in \mathbb{R}^{2 \times 2}$ are boundary layer coefficients with $N_{i,1}^{\mathrm{bl}} > 0$ and $M_{i,1}^{1,\mathrm{bl}} > 0$. It holds $M_{i,2}^{j, \mathrm{bl}} = 0$ for $j=1,2$ for the case of horizontal interfaces considered in this work. For isotropic $\left(\ten{K}_i = k_i \ten{I}\right)$ and orthotropic $\left(\ten{K}_i = \operatorname{diag} (k_{11,i}, \,k_{22,i})\right)$ porous media, we have $N_{i,s}^{\mathrm{bl}} = 0$ and $M_{i,1}^{2,\mathrm{bl}} = 0$.
The boundary layer constants are determined using theory of homogenization and boundary layers and include information on the structure of the porous media in the vicinity of the fracture-matrix interfaces $\gamma_1$ and $\gamma_2$~(Fig.~\ref{fig:fracture-model}). The rigorous derivation of the generalised interface conditions can be found in~\cite[Section~3.2]{Eggenweiler_Rybak_MMS20}. Note that the values of the boundary layer coefficients are dependent on the interface position~\cite[Section~4.1]{Eggenweiler_Rybak_Discacciati_21}. In this work, the interfaces $\gamma_1$ and $\gamma_2$ are located directly below or above the first row of solid inclusions, respectively.

The hybrid-dimensional models with the Beavers--Joseph~\eqref{eq:IC-BJ} and the Beavers--Joseph--Saffman condition~\eqref{eq:IC-BJS} have been developed in~\cite{Rybak-Metzger-20}. Therefore, we refer the reader there for further details. In this work, we derive and analyse an advanced hybrid-dimensional Stokes--Darcy model using the generalised coupling conditions~\eqref{eq:IC-ER}.

\section{Derivation of an advanced hybrid-dimensional model}\label{sec:Reduced-dim-model}
In this section, we provide the main steps for the development of an advanced hybrid-dimensional Stokes--Darcy model. The derivation is based on the vertical averaging approach under the assumptions made in Section~\ref{sec:geometry-assumptions} and follows the same technique as in~\cite{Rybak-Metzger-20} including additional terms due to more sophisticated interface conditions on the top and bottom of the fracture region.

\emph{\textit{Step 1}}: We derive the reduced-dimensional fracture model by averaging the Stokes equations~\eqref{eq:Stokes-macro} in the normal direction to the fracture
\begin{eqnarray*}
     \Omega_\F =  \left\{\vec{x} \in \mathbb{R}^2 \big| \, \vec{x} = \vec{s} +\frac{1}{2} \, \xi \, d \, \vec{n}, \, \vec{s}\in\Gamma, \xi \in [-1,1] \right\},
\end{eqnarray*}
where $d > 0$ denotes the fracture aperture (Fig.~\ref{fig:fracture-model}). Then, we substitute the generalised coupling conditions~\eqref{eq:IC-ER} into the averaged equations.

\emph{\textit{Step 2}}: This intermediate model still contains some terms, where closure relations are needed. 
%We couple the reduced-dimensional fracture model to the porous-medium domains. To obtain the closed model formulation, 
Therefore, we need to make additional assumptions on the flow profile in the fracture, which is typically the case for model derivation by means of vertical averaging, e.g.~\cite{Martin_etal_05, Lesinigo_etal_11, Rybak-Metzger-20}.

\emph{\textit{Step 3}}: We derive transmission conditions to couple the dimensionally reduced fracture model to the porous-medium models.

\emph{\textit{Step 4}}: We formulate the complete hybrid-dimensional model, which contains the Darcy flow equations in the porous-medium domains, the reduced-dimensional fracture model on $\Gamma$ developed in \textit{Step~1} and \textit{Step~2} and the transmission conditions from \textit{Step~3}.

In the following, we define averaged vectorial and scalar quantities 
$$\vec{G}:= \frac{1}{d}\int_{-d/2}^{d/2} \vec{g}_\F \, \d n \quad \text{and} \quad G:= \frac{1}{d}\int_{-d/2}^{d/2} g_\F \, \d n,$$ 
where index $\F$ refers to the fracture domain. A general vector-valued function is expressed in the form
$$\vec{G} = G_{\vec{n}} \vec{n} + G_{\vec{\tau}} \vec{\tau}.$$

\subsection{Step 1 (Vertical averaging)}
Averaging the mass conservation equation~\eqref{eq:Stokes1} across the fracture yields
\begin{eqnarray}
    \vec{v}_\F\cdot \vec{n}| _{\gamma_1} - \vec{v}_\F\cdot \vec{n} | _{\gamma_2} + d \frac{\partial V_\vec{\tau} }{\partial \vec{\tau}} = 0,\label{equ:integratingMass}
\end{eqnarray}
where $V_{\vec{\tau}}$ is the averaged tangential velocity. Substituting the mass conservation equation~\eqref{eq:IC-ER-MC} into~\eqref{equ:integratingMass}, we obtain
\begin{eqnarray}
    \vec{v}_1 \cdot \vec{n} | _{\gamma_1} - \vec{v}_2 \cdot \vec{n} | _{\gamma_2} + d \frac{\partial  V_{\vec{\tau}} }{\partial \vec{\tau}}= 0 \quad \text{on } \Gamma.\label{equ:AveragedMass}
\end{eqnarray}
This equation describes the transport of mass along the fracture.

Integration of the normal component in the momentum conservation equation given in~\eqref{eq:Stokes2} yields
\begin{eqnarray}
         -\left(\mu \frac{\partial \vec{v}_\F}{\partial \vec{n}}\cdot \vec{n}-p_\F\right)\Bigg|_{\gamma_1}  + \left(\mu \frac{\partial \vec{v}_\F}{\partial \vec{n}}\cdot \vec{n} - p_\F\right)\Bigg|_{\gamma_2} 
         = d \left(F_\vec{n} + \mu \frac{\partial^2 V_\vec{n}}{\partial \vec{\tau}^2}\right), \label{equ:averagedMomentumnormaloriginal}
\end{eqnarray}
with the averaged normal velocity 
$V_{\vec{n}}$ and the source term $F_{\vec{n}}$, respectively. 

Substituting the extension of the balance of normal forces~\eqref{eq:IC-ER-BONF} on $\gamma_1$ and $\gamma_2$ into~\eqref{equ:averagedMomentumnormaloriginal}, we obtain 
\begin{eqnarray}
         \bigg( \mu N_{1,s}^{\mathrm{bl}} 
         \frac{\partial \vec{v}_\F}{\partial \vec{n}} \cdot \vec{\tau} +p_1  \bigg)\bigg|_{\gamma_1}    
         -\bigg( \mu N_{2,s}^{\mathrm{bl}} 
         \frac{\partial \vec{v}_\F}{\partial \vec{n}} \cdot \vec{\tau} +p_2  \bigg)\bigg|_{\gamma_2}   = d \left(F_{\vec{n}} + \mu \frac{\partial^2 V_{\vec{n}}}{\partial \vec{\tau}^2}\right).\label{equ:averagedMomentumnormal}
\end{eqnarray}

Analogously, integrating the tangential component of the momentum conservation equation from~\eqref{eq:Stokes2}, we have
\begin{eqnarray}
    - \mu \frac{\partial \vec{v}_\F}{\partial \vec{n}}\cdot \vec{\tau} \bigg|_{\gamma_1} + \mu \frac{\partial \vec{v}_\F}{\partial \vec{n}}\cdot \vec{\tau} \bigg|_{\gamma_2}=  d \left( F_\vec{\tau} +  \mu \frac{\partial^2 V_\vec{\tau}}{\partial \vec{\tau}^2} - \frac{\partial P}{\partial \vec{\tau}} \right),
\label{equ:averagedMomentumtangentialoriginal} \end{eqnarray}
where $P$ is the averaged pressure and $F_{\vec{\tau}}$ is the source term. 

In equations~\eqref{equ:averagedMomentumnormal} and~\eqref{equ:averagedMomentumtangentialoriginal}, additional closure conditions are required to express the terms $\frac{\partial \vec{v}_\F}{\partial \vec{n}} \cdot \vec{\tau}\Big|_{\gamma_i}$ (see derivation in Section~\ref{sec:transmissionconditions}). 

The boundary conditions for the reduced-dimensional fracture model on $\partial \Gamma = \partial \Gamma_D \cup \partial \Gamma_N$ (Fig.~\ref{fig:fracture-model}, right) are determined by averaging the Dirichlet and Neumann boundary conditions~\eqref{eq:BC-Stokes} and are given by
\begin{align}
\label{eq:BC-frac}
   \vec{V} =  \overline{\vec V} \quad \textnormal{on } \partial \Gamma_D, \qquad \frac{1}{d} \int_{-d/2}^{d/2} \ten{T}(\vec{v}_\F, p_\F) {\vec{n}}\; \d n = \overline{\vec{H}} \quad \text{on }\partial \Gamma_N,
\end{align}
where $\overline{\vec V}$ is the averaged Dirichlet boundary value and $\overline{\vec{H}}$ is the averaged Neumann boundary value.

\subsection{Step 2 (Closure relations for dimensionally reduced fracture model)}\label{sec:transmissionconditions}
To get a complete formulation of the dimensionally reduced fracture model given in~\eqref{equ:AveragedMass},~\eqref{equ:averagedMomentumnormal} and~\eqref{equ:averagedMomentumtangentialoriginal}, we need to express $\frac{\partial \vec{v}_\F}{\partial \vec{n}} \cdot \vec{\tau}\Big|_{\gamma_i}$ in~\eqref{equ:averagedMomentumnormal} and~\eqref{equ:averagedMomentumtangentialoriginal} in terms of the averaged quantities $(\vec{V}, P)$ and the interface values $\left(\vec{v}_i, p_i\right)\big|{\gamma_i}$ for $i=1,2$.

Based on~\cite{zimmerman1996hydraulic,rivas2025fluid,Rivas-etal-26}, we consider a quadratic profile for the tangential velocity in the fracture
 \begin{eqnarray}
    \vec{v}_\F \left(n,s\right) \cdot \vec{\tau}= a(s) n^2 + b(s) n +c(s),  \label{eq:tang-quad}
\end{eqnarray}
where the functions $a$, $b$ and $c$ have to be determined. We omit the dependency on $s$, as the expression is valid point-wise. Using~\eqref{eq:tang-quad}, we obtain for the averaged tangential velocity
\begin{eqnarray*}
    V_{\vec{\tau}} = \frac{1}{d}\int^{d/2}_{-d/2} \vec{v}_\F \cdot \vec{\tau} ~\d n = \frac{1}{d}\int^{d/2}_{-d/2} \left( a n^2 + b n + c \right)~\mathrm{d} n = \frac{1}{12} a d^2 + c,  
\end{eqnarray*}
which implies
\begin{eqnarray}
    c = -\frac{1}{12} a d^2 + V_{\vec{\tau}}.\label{eq:c}
\end{eqnarray}
Substituting~\eqref{eq:tang-quad} and~\eqref{eq:c} into the generalisation of the Beavers--Joseph condition~\eqref{eq:IC-ER-BJ} and without loss of generality evaluating the equations at $n= -d/2$ and $n=d/2$ for $i=1,2$, respectively, we get
% at $\gamma_1$ and $\gamma_2$, we get
\begin{subequations}
\label{eq:ab}
\begin{eqnarray}
    \mu \ell  N_{1,1}^\mathrm{bl} (\hspace{+1ex} ad+b\hspace{+1ex}) + \left(\frac{1}{6} ad^2+\frac{1}{2}bd+ V_\vec{\tau}+\frac{\ell^2}{\mu } ( \ten{M}_1^{\mathrm{bl}} \nabla p_1 ) \cdot \vec{\tau} \right) =0, \\
    \mu \ell N_{2,1}^\mathrm{bl} (-ad+b) - \left(\frac{1}{6} ad^2-\frac{1}{2}bd+ V_\vec{\tau}+\frac{\ell^2}{\mu } ( \ten{M}_2^{\mathrm{bl}} \nabla p_2 ) \cdot \vec{\tau} \right) =0.
\end{eqnarray}
\end{subequations}
The solution of the linear system ~\eqref{eq:ab} with respect to $a$ and $b$ yields
\begin{eqnarray}
\label{eq:a}
    a=   \frac{-\left(d+\mu \ell  N_{1,1}^\mathrm{bl}  +\mu \ell N_{2,1}^\mathrm{bl}  \right)V_\vec{\tau}- \left(\mu \ell N_{2,1}^\mathrm{bl}+\frac{1}{2}d\right)\left(\frac{\ell^2}{\mu } ( \ten{M}_1^{\mathrm{bl}} \nabla p_1 ) \cdot \vec{\tau} |_{\gamma_1}\right) - \left(\mu \ell  N_{1,1}^\mathrm{bl} +\frac{1}{2}d\right)\left(\frac{\ell^2}{\mu } ( \ten{M}_2^{\mathrm{bl}} \nabla p_2 ) \cdot \vec{\tau} |_{\gamma_2} \right)}{d\left(\frac{1}{6}d^2 + \frac{2}{3} d\mu \ell  N_{1,1}^\mathrm{bl} +\frac{2}{3}d \mu \ell N_{2,1}^\mathrm{bl}   + 2 \mu^2 \ell^2  N_{1,1}^\mathrm{bl}  N_{2,1}^\mathrm{bl} \right)},
    \\
    \label{eq:b}
    b=\frac{ \left(\mu \ell  N_{1,1}^\mathrm{bl}  -\mu \ell N_{2,1}^\mathrm{bl}  \right)V_\vec{\tau}- \left(\mu \ell N_{2,1}^\mathrm{bl}+\frac{1}{6}d\right)\left(\frac{\ell^2}{\mu } ( \ten{M}_1^{\mathrm{bl}} \nabla p_1 ) \cdot \vec{\tau} |_{\gamma_1}\right) + \left(\mu \ell  N_{1,1}^\mathrm{bl} +\frac{1}{6}d\right)\left(\frac{\ell^2}{\mu } ( \ten{M}_2^{\mathrm{bl}} \nabla p_2 ) \cdot \vec{\tau}|_{\gamma_2} \right) }{\frac{1}{6}d^2 + \frac{2}{3}d \mu \ell  N_{1,1}^\mathrm{bl} +\frac{2}{3}d \mu \ell N_{2,1}^\mathrm{bl}  + 2 \mu^2 \ell^2  N_{1,1}^\mathrm{bl}  N_{2,1}^\mathrm{bl} }.\hspace{+7ex}
\end{eqnarray}
Substituting~\eqref{eq:c},~\eqref{eq:a} and~\eqref{eq:b} into~\eqref{eq:tang-quad} and computing the normal derivatives, we obtain the closure conditions
\begin{eqnarray}
   \left. \frac{\partial \vec{v}_\F}{\partial \vec{n}}\cdot \vec{\tau}\right|_{\gamma_1} =   \frac{ -\left(d+2\mu \ell N_{2,1}^\mathrm{bl}  \right)V_\vec{\tau}- \left(2\mu \ell N_{2,1}^\mathrm{bl}+\frac{2}{3}d\right)\left(\frac{\ell^2}{\mu } ( \ten{M}_1^{\mathrm{bl}} \nabla p_1 ) \cdot \vec{\tau} |_{\gamma_1}\right) - \frac{1}{3}d\left(\frac{\ell^2}{\mu } ( \ten{M}_2^{\mathrm{bl}} \nabla p_2 ) \cdot \vec{\tau} |_{\gamma_2}\right)}{\frac{1}{6}d^2 + \frac{2}{3}d \mu \ell  N_{1,1}^\mathrm{bl} +\frac{2}{3} d\mu \ell N_{2,1}^\mathrm{bl}  + 2 \mu^2 \ell^2  N_{1,1}^\mathrm{bl}  N_{2,1}^\mathrm{bl} },\label{eqn:NormalMomentumClosure1}\\
    \left. \frac{\partial \vec{v}_\F}{\partial\vec{n}}\cdot \vec{\tau}\right|_{\gamma_2} =   \frac{ \left(d+2\mu \ell  N_{1,1}^\mathrm{bl}   \right)V_\vec{\tau}+ \frac{1}{3}d\left(\frac{\ell^2}{\mu } ( \ten{M}_1^{\mathrm{bl}} \nabla p_1 ) \cdot \vec{\tau} |_{\gamma_1}\right) + \left(2\mu \ell  N_{1,1}^\mathrm{bl} +\frac{2}{3}d\right)\left(\frac{\ell^2}{\mu } ( \ten{M}_2^{\mathrm{bl}} \nabla p_2 ) \cdot \vec{\tau} |_{\gamma_2}\right)}{\frac{1}{6}d^2 + \frac{2}{3}d \mu \ell  N_{1,1}^\mathrm{bl} +\frac{2}{3} d\mu \ell N_{2,1}^\mathrm{bl}  + 2 \mu^2 \ell^2  N_{1,1}^\mathrm{bl}  N_{2,1}^\mathrm{bl} }.\label{eqn:NormalMomentumClosure2}\hspace{+2.5ex}
\end{eqnarray}
Inserting~\eqref{eqn:NormalMomentumClosure1} and~\eqref{eqn:NormalMomentumClosure2} into~\eqref{equ:averagedMomentumnormal}, we get
\begin{eqnarray}    
         \left.p_1\right|_{\gamma_1}  - \left.p_2\right|_{\gamma_2}  
         -\frac{  \mu N_{1,s}^{\mathrm{bl}}d + \mu N_{2,s}^{\mathrm{bl}}d+2\mu^2 \ell  N_{1,s}^{\mathrm{bl}} N_{2,1}^\mathrm{bl}  +2\mu^2 \ell   N_{2,s}^{\mathrm{bl}} N_{1,1}^\mathrm{bl} }{\frac{1}{6}d^2 + \frac{2}{3}d \mu \ell  N_{1,1}^\mathrm{bl} +\frac{2}{3} d\mu \ell N_{2,1}^\mathrm{bl}  + 2 \mu^2 \ell^2  N_{1,1}^\mathrm{bl}  N_{2,1}^\mathrm{bl} }V_\vec{\tau}&\nonumber\\ 
         - \frac{ 2\mu \ell^3 N_{1,s}^{\mathrm{bl}}  N_{2,1}^\mathrm{bl}+\frac{2}{3}d \ell^2 N_{1,s}^{\mathrm{bl}}+\frac{1}{3}d \ell^2 N_{2,s}^{\mathrm{bl}}  }{\frac{1}{6}d^2 + \frac{2}{3}d \mu \ell  N_{1,1}^\mathrm{bl} +\frac{2}{3} d\mu \ell N_{2,1}^\mathrm{bl}  + 2 \mu^2 \ell^2  N_{1,1}^\mathrm{bl}  N_{2,1}^\mathrm{bl} }( \ten{M}_1^{\mathrm{bl}} \nabla p_1 ) \cdot \vec{\tau}|_{\gamma_1}&\nonumber\\
          -
          \frac{ 2\mu \ell^3 N_{2,s}^{\mathrm{bl}}  N_{1,1}^\mathrm{bl}+\frac{2}{3}d \ell^2 N_{2,s}^{\mathrm{bl}}+\frac{1}{3}d \ell^2 N_{1,s}^{\mathrm{bl}} }{\frac{1}{6}d^2 + \frac{2}{3}d \mu \ell  N_{1,1}^\mathrm{bl} +\frac{2}{3} d\mu \ell N_{2,1}^\mathrm{bl}  + 2 \mu^2 \ell^2  N_{1,1}^\mathrm{bl}  N_{2,1}^\mathrm{bl} }( \ten{M}_2^{\mathrm{bl}} \nabla p_2 ) \cdot \vec{\tau}|_{\gamma_2}
         & = d \left(F_{\vec{n}} + \mu \frac{\partial^2 V_{\vec{n}}}{\partial \vec{\tau}^2}\right)\quad \textnormal{ on }\Gamma.\label{eqn:averagedMomentNormalFinal}
\end{eqnarray}
Substituting~\eqref{eqn:NormalMomentumClosure1} and~\eqref{eqn:NormalMomentumClosure2} into~\eqref{equ:averagedMomentumtangentialoriginal}, we obtain
\begin{eqnarray}
     \frac{\left(2 \mu d+2\mu^2 \ell N_{1,1}^\mathrm{bl}+2\mu^2 \ell N_{2,1}^\mathrm{bl}  \right)V_\vec{\tau}+ \left(2\mu \ell^3 N_{2,1}^\mathrm{bl}+d\ell^2\right) ( \ten{M}_1^{\mathrm{bl}} \nabla p_1 ) \cdot \vec{\tau} |_{\gamma_1}  + \left(2\mu \ell^3  N_{1,1}^\mathrm{bl} +d\ell^2\right)( \ten{M}_2^{\mathrm{bl}} \nabla p_2 ) \cdot \vec{\tau} |_{\gamma_2}}{\frac{1}{6}d^2 + \frac{2}{3}d \mu \ell  N_{1,1}^\mathrm{bl} +\frac{2}{3} d\mu \ell N_{2,1}^\mathrm{bl}  + 2 \mu^2 \ell^2  N_{1,1}^\mathrm{bl}  N_{2,1}^\mathrm{bl} } \nonumber \hspace{+7ex}\\
     =d \left( F_\vec{\tau} +  \mu \frac{\partial^2 V_\vec{\tau}}{\partial \vec{\tau}^2} - \frac{\partial P}{\partial \vec{\tau}} \right)\quad \textnormal{on }\Gamma.\label{eqn:averagedMomentTangentialFinal}
\end{eqnarray}
Equations~\eqref{eqn:averagedMomentNormalFinal} and~\eqref{eqn:averagedMomentTangentialFinal} describe the transfer of momentum along the lower-dimensional interface $\Gamma$.

\subsection{Step 3 (Transmission conditions)}
To get a closed model formulation, we need to couple the dimensionally reduced fracture model~\eqref{equ:AveragedMass},~\eqref{eqn:averagedMomentNormalFinal} and~\eqref{eqn:averagedMomentTangentialFinal} with the boundary conditions~\eqref{eq:BC-frac} to Darcy's flow equations~\eqref{eq:Darcy-law} and~\eqref{eq:BC-Darcy} in the primal form. Therefore, transmission conditions at the interfaces $\gamma_1$ and $\gamma_2$ are needed. 
For this purpose, pressures $p_1$ and $p_2$ on the interfaces $\gamma_1$ and $\gamma_2$ are required.

We make a priori assumptions on the pressure and velocity profiles in the fracture~$\Omega_\F$, following the investigations from \cite{Lesinigo_etal_11, Rybak-Metzger-20, Ruan_Rybak_AMC}. As the aperture of the fracture is significantly smaller than its length ($d\ll L$), we can neglect the pressure jump and therefore assume a constant pressure profile across the fracture leading to
\begin{eqnarray}
    p_\F |_{\gamma_1} =p_\F |_{\gamma_2} = P. \label{equ:constantpressure} 
\end{eqnarray}

For the normal component of velocity, we consider a quadratic profile across the fracture based on the study in~\cite{Ruan_Rybak_AMC}.
Taking into account the mass conservation~\eqref{eq:IC-ER-MC} on interfaces~$\gamma_1$,~$\gamma_2$ and following similar steps as for the tangential velocity, we obtain
\begin{eqnarray}
      \left.\frac{\partial \vec{v}_\F}{\partial \vec{n}}\cdot\vec{n}\right|_{\gamma_1} = -\frac{1}{d}\left(6V_\vec{n}-4\vec{v}_1 \cdot\vec{n}|_{\gamma_1}-2\vec{v}_2\cdot\vec{n}|_{\gamma_2}\right) ,  \quad
      \left.\frac{\partial \vec{v}_\F}{\partial \vec{n}}\cdot\vec{n}\right|_{\gamma_2}  =\frac{1}{d}\left(6V_\vec{n}-2\vec{v}_1\cdot\vec{n}|_{\gamma_1}-4\vec{v}_2\cdot\vec{n}|_{\gamma_2}\right) .
\label{eqn:closureNormal}
\end{eqnarray}
To reconstruct the pressure term on interface~$\gamma_i$ for the porous medium models, we substitute~\eqref{equ:constantpressure} and~\eqref{eqn:closureNormal} into the interface condition~\eqref{eq:IC-ER-BONF} and obtain the transmission conditions
\begin{align}
    p_1|_{\gamma_1} &= \frac{\mu}{d} \left(6V_\vec{n}-4\vec{v}_1\cdot\vec{n}|_{\gamma_1}-2\vec{v}_2\cdot\vec{n}|_{\gamma_2}\right) + P-\mu N_{1,s}^{\mathrm{bl}} 
          \left.\frac{\partial \vec{v}_\F}{\partial \vec{n}} \cdot \vec{\tau}\right|_{\gamma_1},\label{eq:transmission1}\\
    p_2|_{\gamma_2}&=-\frac{\mu}{d} \left(6V_\vec{n}-2\vec{v}_1\cdot\vec{n}|_{\gamma_1}-4\vec{v}_2\cdot\vec{n}|_{\gamma_2}\right) + P- \mu N_{2,s}^{\mathrm{bl}} 
         \left.\frac{\partial \vec{v}_\F}{\partial \vec{n}} \cdot \vec{\tau}\right|_{\gamma_2}. \label{eq:transmission2}
\end{align}

\subsection{Step 4 (Closed model formulation)}
The complete formulation of the advanced hybrid-dimensional Stokes--Darcy model consists of the Darcy flow equations~\eqref{eq:Darcy-law} in its primal form in the porous-medium domains $\Omega_1$ and $\Omega_2$, the reduced-dimensional model for the fracture on $\Gamma$ given in~\eqref{equ:AveragedMass},\eqref{eqn:averagedMomentNormalFinal} and~\eqref{eqn:averagedMomentTangentialFinal}, the transmission conditions~\eqref{eq:transmission1} and~\eqref{eq:transmission2} and the boundary conditions~\eqref{eq:BC-Darcy} and~\eqref{eq:BC-frac} for the Darcy flow equations and the reduced-dimensional fracture model, respectively.

\section{Analysis of the hybrid-dimensional model}\label{sec:ana}
In this section, we prove existence and uniqueness of the weak solution for the developed hybrid-dimensional Stokes--Darcy system. The proof is conducted for isotropic and orthotropic porous medium configurations. Note that for this case the boundary layer coefficients $N^\mathrm{bl}_{i,s} = 0$ for $i=1,2$ in condition~\eqref{eq:IC-ER-BONF}.

\subsection{Weak formulation for the porous-medium subdomain models}

We consider the porous-medium models in the primal form~\eqref{eq:Darcy-primal} and define the following test function spaces and norms
% $\displaystyle \varphi_\PM \in H_\PM:= \left\{\varphi \in H^1(\Omega_\PM):\, \varphi |_{\Gamma^\D_{\PM}}=0 \right\} $
\begin{equation*}
    \mathcal{H}_i:= \left\{\varphi_i \in H^1(\Omega_i):\, \varphi_i |_{\Gamma_{D,i}}=0 \right\}, \qquad  \| \varphi_i\|_{\mathcal{H}_i}^2: =  \| \varphi_i \|^2_{L^2(\Omega_i)} + \| \nabla \varphi_i\|^2_{L^2(\Omega_i)}\quad \forall \varphi_i \in \mathcal{H}_i, \qquad \text{for } i=1,2. 
\end{equation*}
The boundary data~\eqref{eq:BC-Darcy} are assumed to be $\overline{p}_i \in H^{1/2}\left(\Gamma_{D,i}\right)$ and $\overline{v}_i\in H^{-1/2} \left(\Gamma_{N,i}\right)$. Multiplication of~\eqref{eq:Darcy-primal} by a test function $\varphi_i\in \mathcal{H}_i$  and integration over $\Omega_i$ by parts, yields
\begin{eqnarray}
     \int_{\Omega_i} f_i \varphi_i ~\mathrm{d}\vec{x} &=& \int_{\Omega_i} - \nabla \cdot \left( \frac{\ten{K}_i}{\mu} \nabla p_i\right) \varphi_i~d\vec{x}\nonumber\\
       &=&   \int_{\Omega_i} \left( \frac{\ten{K}_i}{\mu} \nabla p_i\right) \cdot \nabla \varphi_i~\mathrm{d}\vec{x}   + \int_{\gamma_i} (\vec{v}_i \cdot \vec{n}_i) \varphi_i ~\mathrm{d} s 
       +\int_{\Gamma_{N,i}} \overline{v}_i \varphi_i ~\mathrm{d} s\label{equ:weakformulationofDarcy},
\end{eqnarray}
where the boundary conditions \eqref{eq:BC-Darcy} are used. 
From the transmission conditions~\eqref{eq:transmission1} and~\eqref{eq:transmission2}, we obtain for isotropic and orthotropic porous media
% \begin{eqnarray}
%    \vec{v}_1\cdot\vec{n}|_{\gamma_1}  &=& \frac{3}{2} V_\vec{n} + \frac{d}{4\mu}P -\frac{1}{2}\vec{v}_2\cdot\vec{n}|_{\gamma_2}-\frac{d}{4\mu} p_1 \label{eq:transmission1update}\\ %\Ins{- \frac{d N_{1,s}^{\mathrm{bl}}}{4}           \left.\frac{\partial \vec{v}_\F}{\partial \vec{n}} \cdot \vec{\tau}\right|_{\gamma_1}}, \label{eq:transmission1update}\\
%     \vec{v}_2\cdot\vec{n}|_{\gamma_2}&=& \frac{3}{2}V_\vec{n}- \frac{d}{4\mu}P-\frac{1}{2}\vec{v}_1\cdot\vec{n}|_{\gamma_1} \label{eq:transmission2update} % +\frac{d}{4\mu}p_2 \Ins{+ \frac{d N_{2,s}^{\mathrm{bl}}}{4}  \left.\frac{\partial \vec{v}_\F}{\partial \vec{n}} \cdot \vec{\tau}\right|_{\gamma_2}}.\label{eq:transmission2update}
% \end{eqnarray}
% \CHECK{Starting from here, we set $N^\mathrm{bl}_{i,s}=0$ for isotropic/orthotropic porous media.}
% Taking into account Eqs.~\eqref{eq:transmission1update}--\eqref{eq:transmission2update}, we obtain
\begin{eqnarray}
    \vec{v}_1\cdot\vec{n}|_{\gamma_1} &=&  V_\vec{n} + \frac{d}{2\mu}P -\frac{d}{3\mu} p_1|_{\gamma_1} -\frac{d}{6\mu} p_2|_{\gamma_2},\label{eq:transmissionHelp2}\\
     \vec{v}_2\cdot\vec{n}|_{\gamma_2} &=&  V_\vec{n} - \frac{d}{2\mu}P +\frac{d}{6\mu} p_1|_{\gamma_1} +\frac{d}{3\mu} p_2|_{\gamma_2}.
     \label{eq:transmissionHelp3}
      % - \vec{v}_1\cdot\vec{n}|_{\gamma_1}+ \vec{v}_2\cdot\vec{n}|_{\gamma_2}&=& -\frac{d}{\mu}P +\frac{d}{2\mu}p_1 +\frac{d}{2\mu}p_2. \label{eq:transmissionHelp}
\end{eqnarray}

Substituting~\eqref{eq:transmissionHelp2} and~\eqref{eq:transmissionHelp3} into the integral term over the interface~$\gamma_i$ from~\eqref{equ:weakformulationofDarcy} for $i=1, 2$, respectively, we get
\begin{eqnarray}
    \int_{\gamma_1} (\vec{v}_1 \cdot \vec{n}_1) \varphi_1 ~\mathrm{d} s&=&-\int_{\gamma_1} (\vec{v}_1 \cdot \vec{n}) \varphi_1 ~\mathrm{d} s\nonumber \\ 
    &=&-\int_{\gamma_1}  V_\vec{n} \varphi_1 ~\mathrm{d} s -  \frac{d}{2\mu}\int_{\gamma_1} P \varphi_1 ~\mathrm{d} s+  \frac{d}{3\mu} \int_{\gamma_1}  p_1 \varphi_1 ~\mathrm{d} s + \frac{d}{6\mu} \int_{\gamma_1}  p_2 \varphi_1 ~\mathrm{d} s, \label{eq:weakPMclosure1}
\end{eqnarray}
and 
\begin{eqnarray}
    \int_{\gamma_2} (\vec{v}_2 \cdot \vec{n}_2) \varphi_2 ~\mathrm{d} s &=&\int_{\gamma_2}  V_\vec{n} \varphi_2 ~\mathrm{d} s - \frac{d}{2\mu} \int_{\gamma_2} P \varphi_2 ~\mathrm{d} s +  \frac{d}{6\mu} \int_{\gamma_2}  p_1 \varphi_2 ~\mathrm{d} s +  \frac{d}{3\mu} \int_{\gamma_2}  p_2 \varphi_2 ~\mathrm{d} s. \label{eq:weakPMclosure2}
\end{eqnarray}

Considering~\eqref{equ:weakformulationofDarcy},~\eqref{eq:weakPMclosure1} and~\eqref{eq:weakPMclosure2}, the weak formulation for porous-medium subdomain models reads: \\
Find $(p_1, p_2) \in \mathcal{H}_1 \times \mathcal{H}_2 $  such that
\begin{eqnarray}
    \mathcal{A}_\PM  (p_1,p_2; \varphi_1, \varphi_2) +  \mathcal{A}_{\gamma_1,\gamma_2} (p_1,p_2; \varphi_1, \varphi_2)+\mathcal{F}_{\gamma_1, \gamma_2}(\vec{V}, P ; \varphi_1, \varphi_2) = \mathcal{L}_\PM (\varphi_1, \varphi_2), 
    \label{equ:WeakUncoupledPM}
\end{eqnarray}
with the following bilinear operators 
\begin{subequations}
\begin{eqnarray}
    \mathcal{A}_\PM (p_1,p_2; \varphi_1, \varphi_2) &:=& \sum_{i=1,2} \int_{\Omega_i} \left( \frac{\ten{K}_i}{\mu} \nabla p_i\right) \cdot \nabla \varphi_i~\mathrm{d}\vec{x} , \\
    \mathcal{A}_{\gamma_1,\gamma_2} (p_1,p_2; \varphi_1, \varphi_2)&:=& \frac{d}{3\mu} \int_{\gamma_1}  p_1 \varphi_1 ~\mathrm{d} s + \frac{d}{6\mu} \int_{\gamma_1}  p_2 \varphi_1 ~\mathrm{d} s  + \frac{d}{6\mu} \int_{\gamma_2}  p_1 \varphi_2 ~\mathrm{d} s +  \frac{d}{3\mu} \int_{\gamma_2}  p_2 \varphi_2 ~\mathrm{d} s, \\
  \mathcal{F}_{\gamma_1, \gamma_2}(\vec{V}, P ; \varphi_1, \varphi_2)&:=& -\int_{\gamma_1}  V_\vec{n} \varphi_1 ~\mathrm{d} s - \frac{d}{2\mu}\int_{\gamma_1}  P \varphi_1 ~\mathrm{d} s
   +\int_{\gamma_2}  V_\vec{n} \varphi_2 ~\mathrm{d} s - \frac{d}{2\mu}\int_{\gamma_2}  P \varphi_2 ~\mathrm{d} s,
\end{eqnarray}
\end{subequations}
and the linear functional
\begin{eqnarray}
\label{eq:L_pm}
    \mathcal{L}_\PM (\varphi_1, \varphi_2) := \sum_{i=1,2}\int_{\Omega_i} f_i \varphi_i ~\mathrm{d}\vec{x} - \sum_{i=1,2}\int_{\Gamma_{N,i}} \overline{v}_i \varphi_i ~\mathrm{d} s. 
\end{eqnarray}

\subsection{Weak formulation for the reduced-dimensional fracture model}
We now derive the weak formulation of the reduced-dimensional fracture model at the complex interface $\Gamma$. We define the test function space equipped with the norm
\begin{eqnarray*}
    \mathcal{H}_\Gamma : = \left\{ \vec{W}\in \left( H^1\left(\Gamma\right)\right)^2: \vec{W}|_{\partial \Gamma_D} = \vec{0}\right\},\qquad \|\vec{W}\|_{\mathcal{H}_\Gamma}^2:= \|\vec{W} \|_{L^2(\Gamma)}^2+\left\|\frac{\partial \vec{W}}{\partial \vec{\tau}} \right\|_{L^2(\Gamma)}^2 \quad \forall\, \vec{W} \in \mathcal{H}_{\Gamma}.
    % \quad \textnormal{and} \quad \Z := L^2 \left(\Gamma\right)
\end{eqnarray*}
% equipped with the norms 
% \begin{eqnarray*}
%     \|\vec{W}\|_{\mathcal{H}_\Gamma}^2:= \|\vec{W} \|_{L^2(\Gamma)}^2+\left\|\frac{\partial \vec{W}}{\partial \vec{\tau}} \right\|_{L^2(\Gamma)}^2 
%     % \quad \textnormal{and} \quad \| \Psi \|_{\Z} :=\| \Psi \|_{L^2(\Gamma)}. 
% \end{eqnarray*}
We consider the  boundary data~\eqref{eq:BC-frac}, which satisfy the assumption $\overline{\vec{V}}\in \left(H^{1/2}\left(\partial \Gamma_D\right)\right)^2$ and $\overline{\vec{H}}
\in \left(H^{-1/2}(\partial \Gamma_N)\right)^2$.

The weak formulation for the averaged mass conservation equation~\eqref{equ:AveragedMass} is obtained in the standard way 
\begin{eqnarray}
    d \int_\Gamma \frac{\partial V_\vec{\tau}}{\partial \vec{\tau}}  \Psi ~\mathrm{d} s   &=& - \int_\Gamma \left(\vec{v}_1\cdot \vec{n} |_{\gamma_1}\right)\Psi~\mathrm{d} s + \int_\Gamma  \left(\vec{v}_2\cdot \vec{n}|_{\gamma_2}\right) \Psi~\mathrm{d} s \nonumber\\
    &=&- \frac{d}{\mu}\int_\Gamma P\Psi~\mathrm{d} s + \frac{d}{2\mu}\int_\Gamma  \left(p_1 |_{\gamma_1}\right) \Psi~\mathrm{d} s +\frac{d}{2\mu}\int_\Gamma  \left(p_2 |_{\gamma_2}\right) \Psi~\mathrm{d} s \qquad \forall \Psi \in L^2(\Gamma),\label{eqn:weakAveragedMass}
\end{eqnarray}
where conditions~\eqref{eq:transmissionHelp2} and~\eqref{eq:transmissionHelp3} are applied.
Following the standard approach, the weak formulation for the normal component of the averaged momentum conservation equation~\eqref{eqn:averagedMomentNormalFinal} for the case of isotropic and orthotropic porous media is obtained
\begin{eqnarray}
        % \int_\Gamma ( p_1 |_{\gamma_1} ) W_\vec{n}  ~\mathrm{d} s 
        %  - \int_\Gamma(p_2  |_{\gamma_2})  W_\vec{n} ~\mathrm{d} s -d\mu\int_\Gamma  \left(\frac{\partial^2 V_{\vec{n}}}{\partial \vec{\tau}^2}\right)W_\vec{n} ~\mathrm{d} s& =& d \int_\Gamma F_{\vec{n}}W_\vec{n} ~\mathrm{d} s \nonumber \\
          \int_\Gamma ( p_1 |_{\gamma_1} ) W_\vec{n}  ~\mathrm{d} s 
         - \int_\Gamma(p_2  |_{\gamma_2})  W_\vec{n} ~\mathrm{d} s+ d \mu\int_\Gamma  \frac{\partial V_\vec{n}}{\partial \vec{\tau}} \frac{\partial W_{\vec{n}}}{\partial \vec\tau }~\mathrm{d} s &=& d \int_\Gamma F_{\vec{n}}W_\vec{n} ~\mathrm{d} s + d\left[ \overline{H}_{\vec n} W_{\vec n}\right]_{\partial \Gamma_N}.
\end{eqnarray}
Analogously, we derive the weak formulation of the tangential component~\eqref{eqn:averagedMomentTangentialFinal}:
\begin{eqnarray}
     d \int_\Gamma \left(\mu \frac{\partial V_\vec{\tau}}{\partial \vec{\tau}} - P \right) \frac{\partial W_{\vec{\tau}}}{\partial \vec{\tau}}~\mathrm{d} s+\frac{\left(2 \mu d+2\mu^2 \ell N_{1,1}^\mathrm{bl}+2\mu^2 \ell N_{2,1}^\mathrm{bl}  \right)}{\frac{1}{6}d^2 + \frac{2}{3}d \mu \ell  N_{1,1}^\mathrm{bl} +\frac{2}{3} d\mu \ell N_{2,1}^\mathrm{bl}  + 2 \mu^2 \ell^2  N_{1,1}^\mathrm{bl}  N_{2,1}^\mathrm{bl} } \int_\Gamma V_{\vec{\tau}}  W_{\vec{\tau}}~\mathrm{d} s\nonumber \hspace{+10ex}\\
     +\frac{ \left(2\mu \ell^3 N_{2,1}^\mathrm{bl}+d\ell^2\right)  }{\frac{1}{6}d^2 + \frac{2}{3}d \mu \ell  N_{1,1}^\mathrm{bl} +\frac{2}{3} d\mu \ell N_{2,1}^\mathrm{bl}  + 2 \mu^2 \ell^2  N_{1,1}^\mathrm{bl}  N_{2,1}^\mathrm{bl} }\int_\Gamma \left(\left.( \ten{M}_1^{\mathrm{bl}} \nabla p_1 ) \cdot \vec{\tau}\right|_{\gamma_1}\right) W_{\vec{\tau}}~\mathrm{d} s \nonumber \\
     +\frac{ \left(2\mu \ell^3  N_{1,1}^\mathrm{bl} +d\ell^2\right) }{\frac{1}{6}d^2 + \frac{2}{3}d \mu \ell  N_{1,1}^\mathrm{bl} +\frac{2}{3} d\mu \ell N_{2,1}^\mathrm{bl}  + 2 \mu^2 \ell^2  N_{1,1}^\mathrm{bl}  N_{2,1}^\mathrm{bl} } \int_\Gamma \left(\left.( \ten{M}_2^{\mathrm{bl}} \nabla p_2 ) \cdot \vec{\tau}\right|_{\gamma_2}\right) W_{\vec{\tau}}~\mathrm{d} s \nonumber\\
     =d \int_\Gamma  F_\vec{\tau} W_{\vec{\tau}}~\mathrm{d} s + d\left[ \overline{H}_{\vec \tau} W_{\vec \tau}\right]_{\partial \Gamma_N}\, .
\end{eqnarray}
The weak formulation of the reduced-dimensional fracture model is given by: \\
Find $(\vec{W}, P) \in \mathcal{H}_\Gamma\times \Z $  such that
\begin{subequations}
\begin{align}
        \mathcal{A}_\Gamma\left( \vec{V}; \vec{W}\right)   + \mathcal{B}_\Gamma\left(\vec{W}; P\right)  &= \mathcal{F}_\Gamma( p_1, p_2;\vec{W})+\mathcal{L}_\Gamma\left( \vec{W}\right) &&\forall\vec{W} \in \mathcal{H}_\Gamma,  \label{equ:weakUncoupledgamma1} \\
        %next
    \mathcal{B}_\Gamma \left( \vec{V}; \Psi \right) - \mathcal{E}_\Gamma\left(P;\Psi \right) &= \mathcal{G}_\Gamma (p_1, p_2; \Psi ) &&\forall \Psi\hspace{+0.5ex}\in\hspace{+0.2ex} \Z, \label{equ:weakUncoupledgamma2}
\end{align}\label{equ:weakUncoupledgamma}
\end{subequations}
where the corresponding bilinear operators and linear functional are defined as 
\begin{subequations}
\begin{eqnarray}
    \mathcal{A}_\Gamma(\vec{V};\vec{W})&:=&  d\mu\int_\Gamma  \frac{\partial V_\vec{n}}{\partial \vec{\tau}} \frac{\partial W_{\vec{n}}}{\partial \vec\tau }~\mathrm{d} s  +d\mu \int_\Gamma   \frac{\partial V_\vec{\tau}}{\partial \vec{\tau}}\frac{\partial W_{\vec{\tau}}}{\partial \vec{\tau}}~\mathrm{d} s \nonumber\\
    &\quad&+\frac{\left(2 \mu d+2\mu^2 \ell N_{1,1}^\mathrm{bl}+2\mu^2 \ell N_{2,1}^\mathrm{bl}  \right)}{\frac{1}{6}d^2 + \frac{2}{3}d \mu \ell  N_{1,1}^\mathrm{bl} +\frac{2}{3} d\mu \ell N_{2,1}^\mathrm{bl}  + 2 \mu^2 \ell^2  N_{1,1}^\mathrm{bl}  N_{2,1}^\mathrm{bl} } \int_\Gamma V_{\vec{\tau}}  W_{\vec{\tau}}~\mathrm{d} s,\\ 
    \mathcal{B}_\Gamma(\vec{W};P) &:=&- \int_\Gamma   d P \frac{\partial W_{\vec{\tau}}}{\partial \vec{\tau}}~\mathrm{d} s, \\
    \mathcal{E}_\Gamma(P, \Psi) &:=&\int_\Gamma \frac{d}{\mu} P  \Psi ~\mathrm{d} s,\\
     \mathcal{F}_\Gamma(p_1,p_2; \vec{W}) &:=& - \int_\Gamma ( p_1 |_{\gamma_1} ) W_\vec{n}  ~\mathrm{d} s 
         + \int_\Gamma(p_2  |_{\gamma_2})  W_\vec{n} ~\mathrm{d} s
     \nonumber\\
     &\quad&-  \frac{ \left(2\mu \ell^3 N_{2,1}^\mathrm{bl}+d\ell^2\right)  }{\frac{1}{6}d^2 + \frac{2}{3}d \mu \ell  N_{1,1}^\mathrm{bl} +\frac{2}{3} d\mu \ell N_{2,1}^\mathrm{bl}  + 2 \mu^2 \ell^2  N_{1,1}^\mathrm{bl}  N_{2,1}^\mathrm{bl} }\int_\Gamma \left(\left.( \ten{M}_1^{\mathrm{bl}} \nabla p_1 ) \cdot \vec{\tau}\right|_{\gamma_1}\right) W_{\vec{\tau}}~\mathrm{d} s \nonumber \hspace{+5ex}\\
     &\quad&-\frac{ \left(2\mu \ell^3  N_{1,1}^\mathrm{bl} +d\ell^2\right) }{\frac{1}{6}d^2 + \frac{2}{3}d \mu \ell  N_{1,1}^\mathrm{bl} +\frac{2}{3} d\mu \ell N_{2,1}^\mathrm{bl}  + 2 \mu^2 \ell^2  N_{1,1}^\mathrm{bl}  N_{2,1}^\mathrm{bl} } \int_\Gamma \left(\left.( \ten{M}_2^{\mathrm{bl}} \nabla p_2 ) \cdot \vec{\tau}\right|_{\gamma_2}\right) W_{\vec{\tau}}~\mathrm{d} s, \\
    \mathcal{G}_\Gamma(p_1,p_2;  \Psi) &:=& -\frac{d}{2\mu}\int_\Gamma  \left(p_1 |_{\gamma_1}\right) \Psi~\mathrm{d} s -\frac{d}{2\mu}\int_\Gamma  \left(p_2 |_{\gamma_2}\right) \Psi~\mathrm{d} s,\\
     \mathcal{L}_\Gamma\left( \vec{W}\right)&:=&d \int_\Gamma  F_\vec{n}W_{\vec{n}}~\mathrm{d} s +d\int_\Gamma    F_\vec{\tau} W_{\vec{\tau}}~\mathrm{d} s+ d\left[ \overline{H}_{\vec n} W_{\vec n}\right]_{\partial \Gamma_N} + d\left[ \overline{H}_{\vec \tau} W_{\vec \tau}\right]_{\partial \Gamma_N}. \label{eq:L_Gamma}
\end{eqnarray}
\end{subequations}

\subsection{Weak formulation for the hybrid-dimensional Stokes--Darcy model}

In this section, we provide the weak formulation for the hybrid-dimensional problem.
% and then prove the well-posedness by following the approach in~\cite{Lesinigo_etal_11, Rybak-Metzger-20, Ruan_Rybak_AMC}.  
We define the test function space
\begin{eqnarray*}
    \vec \theta : =( \varphi_1, \varphi_2, \vec{W}, \Psi) \in \mathcal{H}_\HYB := \mathcal{H}_1 \times \mathcal{H}_2 \times \mathcal{H}_\Gamma\times\Z
\end{eqnarray*}
equipped with the norm
\begin{eqnarray*}
    \| \vec \theta  \|_{\mathcal{H}_\HYB}^2 :=  \| \varphi_1 \|_{\mathcal{H}_1}^2 + \| \varphi_2 \|_{\mathcal{H}_2}^2 + \| \vec{W}\|_{\mathcal{H}_\Gamma}^2 + \|\Psi\|_{\Z}^2.
\end{eqnarray*}
Taking into account equations~\eqref{equ:WeakUncoupledPM} and~\eqref{equ:weakUncoupledgamma}, we obtain the weak formulation of the coupled hybrid-dimensional model:\\
Find $\vec{\zeta}:=(p_1, p_2, \vec{V},P)\in \mathcal{H}_\HYB $  such that
\begin{eqnarray}
    \mathcal{A}_\HYB (\vec{\zeta}; \vec{\theta}  )=  \mathcal{L}_\HYB(\vec{\theta}), \quad \forall \vec{\theta} \in \mathcal{H}_\HYB, \label{eq:weakformCoupledHybrid}
\end{eqnarray}
where the following bilinear operator and linear functional are defined
\begin{subequations}
\begin{eqnarray}
\label{eq:weakformHybrid1}
      \mathcal{A}_\HYB (\vec{\zeta}; \vec{\theta}  )&:=&  \mathcal{A}_\Gamma\left( \vec{V}; \vec{W}\right)   + \mathcal{B}_\Gamma\left(\vec{W}; P\right)  -  \mathcal{F}_\Gamma( p_1, p_2;\vec{W})
       -\mathcal{B}_\Gamma \left( \vec{V}; \Psi \right) + \mathcal{E}_\Gamma\left(P;\Psi \right)+ \mathcal{G}_\Gamma (p_1, p_2; \Psi ) \nonumber
      \\
      &\quad&
      +\mathcal{A}_\PM  (p_1,p_2; \varphi_1, \varphi_2) +  \mathcal{A}_{\gamma_1,\gamma_2} (p_1,p_2; \varphi_1, \varphi_2)+\mathcal{F}_{\gamma_1, \gamma_2}(\vec{V}, P ; \varphi_1, \varphi_2)  ,\\
    \mathcal{L}_\HYB(\vec{\theta})&:=& \mathcal{L}_\Gamma\left( \vec{W}\right) +\mathcal{L}_\PM (\varphi_1, \varphi_2). \label{eq:weakformHybrid2}
\end{eqnarray}
\end{subequations}

\begin{remark}
    We handle non-homogeneous Dirichlet boundary conditions $(\overline{p}_1, \overline{p}_2, \overline{\vec V}) \in H^{1/2}\left(\Gamma_{D,1}\right)\times H^{1/2}(\Gamma_{D,2}) \times \left(H^{1/2}\left(\partial\Gamma_{D}\right)\right)^2$ by introducing a corresponding lifting $\Xi \in H^1(\Omega_1) \times H^1(\Omega_2) \times \left(H^1(\Gamma)\right)^2$ of $(\hat{p}_1, \hat{p}_2, \hat{\vec{V}})\in H^1(\Omega_1) \times H^1(\Omega_2) \times \left(H^1(\Gamma)\right)^2$ such that $((\hat{p}_1, \hat{p}_2, \hat{\vec{V}})-\Xi,P)=\vec{\zeta}\in \mathcal{H}_\HYB$ in the formulation~\eqref{eq:weakformCoupledHybrid}.
\end{remark}

\subsection{Well-posedness analysis}
In this section, first we introduce the auxiliary results needed to prove the well-posedness of the hybrid-dimensional model developed in this work and then provide the proof of Theorem~\ref{theo:well-posedness}. Due to complexity of the hybrid-dimensional model, we introduce the following notations for better readability
\begin{subequations}
\label{eq:abbr}
\begin{align}
    B(d) &:= \frac{\left(2 \mu d+2\mu^2 \ell N_{1,1}^\mathrm{bl}+2\mu^2 \ell N_{2,1}^\mathrm{bl}  \right)}{\frac{1}{6}d^2 + \frac{2}{3}d \mu \ell  N_{1,1}^\mathrm{bl} + \frac{2}{3}d\mu \ell N_{2,1}^\mathrm{bl}  + 2 \mu^2 \ell^2  N_{1,1}^\mathrm{bl}  N_{2,1}^\mathrm{bl} }, \label{eq:abbr-B}\\
    b_1(d) &:= \frac{ \left(2\mu \ell^3 N_{2,1}^\mathrm{bl}+d\ell^2\right)\MI}{\frac{1}{6}d^2 + \frac{2}{3}d \mu \ell  N_{1,1}^\mathrm{bl} +\frac{2}{3} d\mu \ell N_{2,1}^\mathrm{bl}  + 2 \mu^2 \ell^2  N_{1,1}^\mathrm{bl}  N_{2,1}^\mathrm{bl} }, \label{eq:abbr-b1}\\
    b_2(d) &:= \frac{ \left(2\mu \ell^3  N_{1,1}^\mathrm{bl} +d\ell^2\right)\MII }{\frac{1}{6}d^2 + \frac{2}{3}d \mu \ell  N_{1,1}^\mathrm{bl} +\frac{2}{3} d\mu \ell N_{2,1}^\mathrm{bl}  + 2 \mu^2 \ell^2  N_{1,1}^\mathrm{bl}  N_{2,1}^\mathrm{bl} }, \label{eq:abbr-b2}\\
    C_{T_2}(d)
    &:= \frac{1}{C_I}
    \min\left\{
        B(d),\,
        \sqrt{2d\mu B(d)}
    \right\}, \label{eq:abbr-CT2}
    \\
    C_{T_{3,4}}(d)
    &:= \sum_{i\in\{1,2\}}
    \frac{
        b_i(d)^2\mu C_{\vec{\tau},i}^2
        \tilde{C}_{P,i}
    }{
        k_{\min,i}
    }. \label{eq:abbr-CT34}
\end{align}
\end{subequations}

\begin{theorem}[Well-posedness of the hybrid-dimensional model]
\label{theo:well-posedness}
    If the inequality
\begin{equation}
    \frac{C_{T_{3,4}}(d)}{4C_{T_2}(d)} < 1
    \label{eq:coercivity-condition}
\end{equation}
is satisfied, the solution of the hybrid-dimensional fracture model~\eqref{eq:weakformCoupledHybrid} exists and is unique.
\end{theorem}

\subsubsection{Auxiliary results}
Here, we summarise the estimates, which are used in the proof of Theorem~\ref{theo:well-posedness}.
% the trace theorem~\cite[Theorem~I.1.5]{girault2012finite}, the trace theorem for the tangential derivative~\cite[Equation (3.18)]{Eggenweiler_Rybak_Discacciati_21}, the Sobolev embedding  and the Poincar{\'e} inequality~\cite[Theorem~I.1.1]{girault2012finite}.
\begin{lemma}
[Trace theorem \cite{girault2012finite}, Theorem~I.1.5]\label{thm:tracetheorem}
Let $\Omega \subset \mathbb{R}^n$ be a bounded Lipschitz domain with boundary $\partial \Omega$. Then the mapping $f\to f|_{\partial \Omega}$  has a unique linear continuous extension as an operator from 
% \begin{eqnarray*}
    $H^1(\Omega)$ onto $H^{1/2}(\partial \Omega)$ implying
% \end{eqnarray*}
% and let $p\ge 1$ and $s\ge0$ be two real numbers such that $s\le k+1$, $s-\frac{1}{p}=l+\sigma$, where $l\ge0$ is an integer and $\sigma$.
\begin{eqnarray}
        \exists C_{0}>0 &\textnormal{s.t.}& \| f|_{\partial \Omega} \|_{L^2(\partial \Omega)} \le \|f |_{\partial \Omega}\|_{H^{1 / 2}(\partial \Omega)} \le C_{0} \|f \|_{H^1 (\Omega)} \qquad \forall f \in H^1(\Omega).\label{equ:tracegeneral}
\end{eqnarray}
\end{lemma}

\begin{lemma}[Trace theorem for the tangential derivative \cite{Eggenweiler_Rybak_Discacciati_21}, Equation (3.18)]
\label{thm:tracetheoremtangential}
Let $\Omega\subset\mathbb{R}^2$ be a bounded Lipschitz domain with boundary $\partial \Omega$.
For $f\in H^1(\Omega)$, the tangential derivative of its trace
belongs to $H^{-1/2}(\partial\Omega)$ and satisfies
\begin{eqnarray}
        \exists C_{\vec{\tau}}>0 &\textnormal{s.t.}& \| (\nabla f) \cdot\vec{\tau}|_{\partial \Omega} \|_{H^{-1 / 2}(\partial \Omega)} \le C_{\vec{\tau}} \|f \|_{H^1 (\Omega)} \qquad \forall f \in H^1(\Omega).\label{equ:tracetangential}
\end{eqnarray}
% Consequently, for every $w\in H^{1/2}(\partial\Omega)$,
% \begin{equation}
%     \left|
%     \left\langle
%     (\nabla f) \cdot\vec{\tau}|_{\partial \Omega},w
%     \right\rangle_
%     {H^{-1/2}(\partial\Omega),H^{1/2}(\partial\Omega)}
%     \right|
%     \le
%     C_{\boldsymbol{\tau}}
%     \|f\|_{H^1(\Omega)}
%     \|w\|_{H^{1/2}(\partial\Omega)}.
% \end{equation}
\end{lemma}

% \begin{lemma}
% [Trace theorem for the tangential derivative]\label{thm:tracetheoremtangential}
% Let $\Omega \subset \mathbb{R}^2$ be a bounded Lipschitz domain with boundary~$\partial \Omega$. Then the following statement is valid
% \begin{eqnarray}
%         \exists C_{\vec{\tau}}>0 &\textnormal{s.t.}& \| (\nabla f) \cdot\vec{\tau}|_{\partial \Omega} \|_{H^{-1 / 2}(\partial \Omega)} \le C_{\vec{\tau}} \|f \|_{H^1 (\Omega)}
%          \qquad \forall f \in H^1(\Omega).\label{equ:tracetangential}
% \end{eqnarray}
% \end{lemma}
% \Dis{Check the connection between $L^2$ and $H^{-1/2}$}

\begin{lemma}[Sobolev embedding on the interface \cite{McLean-01}, page 76]
\label{lem:embedding-interface}
Let $\Gamma\subset\mathbb{R}^2$ be a bounded one-dimensional
Lipschitz manifold. Then the embedding
\[
H^1(\Gamma)\hookrightarrow H^{1/2}(\Gamma)
\]
is continuous, i.e.
\begin{eqnarray}
    \exists \Ce>0 &\textnormal{s.t.}&\|w\|_{H^{1/2}(\Gamma)}
    \le
    \Ce\|w\|_{H^1(\Gamma)}
    \qquad
    \forall w\in H^1(\Gamma).
\end{eqnarray}
\end{lemma}

\begin{lemma}[Sobolev interpolation inequality \cite{VanSchaftingen-23}, Equation (2)]
\label{lem:Sobolev-interpolation}
    Let \(\Gamma\subset\mathbb{R}^n\) be a bounded Lipschitz domain. Then
    \begin{eqnarray}
        \exists C_I > 0 &\textnormal{s.t.}& \|w\|^2_{H^{1/2}(\Gamma) } \leq C_I \|w\|_{L^2(\Gamma)}\|w\|_{H^1(\Gamma)} \qquad \forall w \in H^1(\Gamma).
    \end{eqnarray}
\end{lemma}

\begin{lemma}[Poincar{\'e} inequality \cite{boffi2013mixed}, Section~1.2] \label{thm:poicare} Let the domain $\Omega$ be connected and bounded. Then 
\begin{eqnarray}
\exists C_P > 0 &\textnormal{s.t.}&
    \| f \|_{L^2(\Omega)} \le  C_{P} \| \nabla f \|_{L^2(\Omega)} \qquad \forall f \in \{H^1(\Omega) : f|_{\Gamma_D} = 0\},
\end{eqnarray}
where $\Gamma_D \neq \emptyset$.
\end{lemma}

Taking $\|f\|_{H^1(\Omega)} ^2 = \|f\|_{L^2(\Omega)} ^2 +\|\nabla f\|_{L^2(\Omega)} ^2 $ into account, we get
\begin{eqnarray}
    \|f\|_{H^1(\Omega)} ^2 \le \tilde{C}_{P} \|\nabla f\|_{L^2(\Omega)} ^2, \quad \tilde{C}_{P}=1+C_{P}^2>1. \label{equ:auxiliarypoincare}
\end{eqnarray}

\subsubsection{Proof of Theorem~\ref{theo:well-posedness}}
To prove well-posedness of problem~\eqref{eq:weakformCoupledHybrid}, we verify the conditions in the Lax--Milgram theorem, e.g.~\cite[Theorem~4.1.6]{boffi2013mixed}. 
The continuity of the linear operator $\mathcal{L}_\HYB$ is straightforward. Now, we verify the continuity of the bilinear form $\mathcal{A}_\HYB$:
\begin{align}
    |\mathcal{A}_\HYB (\vec{\zeta}; \vec{\theta}  )|&\leq  |\mathcal{A}_\Gamma\left( \vec{V}; \vec{W}\right)|   + |\mathcal{B}_\Gamma\left(\vec{W}; P\right)|  +  |\mathcal{F}_\Gamma( p_1, p_2;\vec{W})|
       +|\mathcal{B}_\Gamma \left( \vec{V}; \Psi \right)| + |\mathcal{E}_\Gamma\left(P;\Psi \right)| + |\mathcal{G}_\Gamma (p_1, p_2; \Psi )| \nonumber
      \\
      &\quad
      +|\mathcal{A}_\PM  (p_1,p_2; \varphi_1, \varphi_2)| +  |\mathcal{A}_{\gamma_1,\gamma_2} (p_1,p_2; \varphi_1, \varphi_2)| +| \mathcal{F}_{\gamma_1, \gamma_2}(\vec{V}, P ; \varphi_1, \varphi_2)|\label{eq:continuity}
      % &\leq C_{\mathcal{A}_\HYB} \|\vec{\zeta}\|_{\mathcal{H}_\HYB} \|\vec{\theta}\|_{\mathcal{H}_\HYB}
\end{align}
for all $\vec{\zeta}= (p_1, p_2, \vec{V}, P), \;\vec{\theta} = (\varphi_1, \varphi_2, \vec{W}, \Psi) \in \mathcal{H}_\HYB$.
Below, we bound the terms appearing in the right-hand side of~\eqref{eq:continuity}.
Applying the Cauchy--Schwarz inequality, the ellipticity conditions in~\eqref{equ:boundedK}, and the trace theorem (Lemma~\ref{thm:tracetheorem}), we obtain
\begin{eqnarray}
    &\quad&|\mathcal{A}_\PM  (p_1,p_2; \varphi_1, \varphi_2)| +  |\mathcal{A}_{\gamma_1,\gamma_2} (p_1,p_2; \varphi_1, \varphi_2) | \\
    &\le& \sum_{i=1,2} \frac{k_{\max,i}}{\mu} \left\|\nabla p_i\right\|_{L^2(\Omega_i)}\left\|\nabla \varphi_i\right\|_{L^2(\Omega_i)} \nonumber \\
    &\quad&+\sum_{i=1,2} \frac{d}{3\mu} \left\| p_i\right\|_{L^2(\Gamma)}\left\| \varphi_i\right\|_{L^2(\Gamma)}
    + \frac{d}{6\mu} \left\| p_2\right\|_{L^2(\Gamma)}\left\| \varphi_1\right\|_{L^2(\Gamma)}
    + \frac{d}{6\mu} \left\| p_1\right\|_{L^2(\Gamma)}\left\| \varphi_2\right\|_{L^2(\Gamma)}\nonumber\\
    &\le& \sum_{i=1,2} \frac{k_{\max,i}}{\mu} \left\| p_i\right\|_{\mathcal{H}_i}\left\|\varphi_i\right\|_{\mathcal{H}_i} 
    +\sum_{i=1,2} \frac{d C_{0,i}^2}{3\mu} \left\| p_i\right\|_{\mathcal{H}_i}\left\| \varphi_i\right\|_{\mathcal{H}_i}\nonumber \\
    &\quad&+ \frac{dC_{0,1}C_{0,2}}{6\mu} \left\| p_2\right\|_{\mathcal{H}_2}\left\| \varphi_1\right\|_{\mathcal{H}_1}
    + \frac{dC_{0,1}C_{0,2}}{6\mu} \left\| p_1\right\|_{\mathcal{H}_1}\left\| \varphi_2\right\|_{\mathcal{H}_2}\nonumber\\
    &\le& C_{\mathcal{A}_\PM}\left(\left\| p_1\right\|_{\mathcal{H}_1}+\left\| p_2\right\|_{\mathcal{H}_2}\right)\left(\left\| \varphi_1\right\|_{\mathcal{H}_1}+\left\| \varphi_2\right\|_{\mathcal{H}_2}\right), \label{eq:continuity1}
\end{eqnarray}
with $\displaystyle
    C_{\mathcal{A}_\PM}:= \max\left\{\frac{k_{\max,1}}{\mu}  + \frac{d C_{0,1}^2}{3\mu},\; \frac{k_{\max,2}}{\mu}  + \frac{d C_{0,2}^2}{3\mu},\;\frac{dC_{0,1}C_{0,2}}{6\mu} \right\} .
$

Applying the Cauchy--Schwarz inequality to the first term in the right-hand side of equation~\eqref{eq:continuity}, we get
\begin{eqnarray}
\label{eq:continuity2}
    |\mathcal{A}_\Gamma(\vec{V};\vec{W})|&\le& d\mu \left\| \frac{\partial V_\vec{n}}{\partial \vec{\tau}}\right\|_{L^2(\Gamma)} \left\|\frac{\partial W_{\vec{n}}}{\partial \vec\tau }\right\|_{L^2(\Gamma)}  +d\mu \left\| \frac{\partial V_\vec{\tau}}{\partial \vec{\tau}}\right\|_{L^2(\Gamma)} \left\|\frac{\partial W_{\vec{\tau}}}{\partial \vec\tau }\right\|_{L^2(\Gamma)}  \nonumber \\
    &\quad &+\frac{\left(2 \mu d+2\mu^2 \ell N_{1,1}^\mathrm{bl}+2\mu^2 \ell N_{2,1}^\mathrm{bl}  \right)}{\frac{1}{6}d^2 + \frac{2}{3}d \mu \ell  N_{1,1}^\mathrm{bl} +\frac{2}{3} d\mu \ell N_{2,1}^\mathrm{bl}  + 2 \mu^2 \ell^2  N_{1,1}^\mathrm{bl}  N_{2,1}^\mathrm{bl} }\left\|  V_\vec{\tau}\right\|_{L^2(\Gamma)} \left\|W_{\vec{\tau}}\right\|_{L^2(\Gamma)} \nonumber\\
    &\le& C_{\mathcal{A}_\Gamma} \|\vec{V}\|_{\mathcal{H}_\Gamma}\|\vec{W}\|_{\mathcal{H}_\Gamma},
\end{eqnarray}
with 
$ \displaystyle
  C_{\mathcal{A}_\Gamma}: = d\mu +\frac{\left(2 \mu d+2\mu^2 \ell N_{1,1}^\mathrm{bl}+2\mu^2 \ell N_{2,1}^\mathrm{bl}  \right)}{\frac{1}{6}d^2 + \frac{2}{3}d \mu \ell  N_{1,1}^\mathrm{bl} +\frac{2}{3} d\mu \ell N_{2,1}^\mathrm{bl}  + 2 \mu^2 \ell^2  N_{1,1}^\mathrm{bl}  N_{2,1}^\mathrm{bl}}. 
$

It is evident that the following inequalities hold
\begin{eqnarray}
\label{eq:continuity3}
     |\mathcal{B}_\Gamma(\vec{W};P)| \le d \|P\|_{\Z} \|\vec{W}\|_{\mathcal{H}_\Gamma}, \quad |\mathcal{E}_\Gamma(P, \Psi)| \le \frac{d}{\mu} \|P\|_{\Z} \|\Psi\|_{\Z}.
\end{eqnarray}
Next, we bound the last term in equation~\eqref{eq:continuity}. Taking into account the Cauchy--Schwarz inequality and the trace theorem (Lemma~\ref{thm:tracetheorem}), we obtain
\begin{eqnarray}
\label{eq:continuity4}
     |\mathcal{F}_{\gamma_1, \gamma_2}(\vec{V}, P ; \varphi_1, \varphi_2)|&\le& \left\|  V_\vec{n}\right\|_{L^2(\Gamma)} \left\|\varphi_1\right\|_{L^2(\Gamma)} +\left\|  V_\vec{n}\right\|_{L^2(\Gamma)} \left\|\varphi_2\right\|_{L^2(\Gamma)} +
    \frac{d}{2\mu} \left\|  P\right\|_{L^2(\Gamma)} \left\|\varphi_1\right\|_{L^2(\Gamma)} \nonumber\\ 
    &\quad&+\frac{d}{2\mu} \left\|  P\right\|_{L^2(\Gamma)} \left\|\varphi_2\right\|_{L^2(\Gamma)} \nonumber\\
    &\le&C_{0,1}\left\|  \vec{V}\right\|_{\mathcal{H}_\Gamma} \left\|\varphi_1\right\|_{\mathcal{H}_1} +C_{0,2}\left\|  \vec{V}\right\|_{\mathcal{H}_\Gamma}  \left\|\varphi_2\right\|_{\mathcal{H}_2}  +
    \frac{dC_{0,1}}{2\mu} \left\|  P\right\|_{\Z} \left\|\varphi_1\right\|_{\mathcal{H}_1}\nonumber\\ 
    &\quad&+\frac{dC_{0,2}}{2\mu} \left\|  P\right\|_{\Z} \left\|\varphi_2\right\|_{\mathcal{H}_2}\nonumber\\
    &\le& C_{\mathcal{F}_{\gamma_1,\gamma_2}}\left(\left\|  \vec{V}\right\|_{\mathcal{H}_\Gamma}+\left\|  P\right\|_{\Z}\right)\left(\left\|\varphi_1\right\|_{\mathcal{H}_1}+\left\|\varphi_2\right\|_{\mathcal{H}_2}\right),
\end{eqnarray}
where $\displaystyle C_{\mathcal{F}_{\gamma_1,\gamma_2}}:=\max\left\{ C_{0,1},\;C_{0,2}\;\frac{dC_{0,1}}{2\mu},\;\frac{dC_{0,2}}{2\mu}\right\}$.

In a similar manner, we estimate
\begin{eqnarray}
\label{eq:continuity5}
     |\mathcal{G}_\Gamma(p_1,p_2;  \Psi)| &\le& \frac{d}{2\mu}\left\| p_1\right\|_{L^2(\Gamma)}\left\| \Psi\right\|_{L^2(\Gamma)}+\frac{d}{2\mu}\left\| p_2\right\|_{L^2(\Gamma)}\left\| \Psi\right\|_{L^2(\Gamma)}\nonumber\\
     &\le& \frac{dC_{0,1}}{2\mu}\left\| p_1\right\|_{\mathcal{H}_1}\left\| \Psi\right\|_{\Z}+\frac{dC_{0,2}}{2\mu}\left\| p_2\right\|_{\mathcal{H}_2}\left\| \Psi\right\|_{\Z}\nonumber\\
     &\le& C_{\mathcal{G}_\Gamma}\left(\left\| p_1\right\|_{\mathcal{H}_1}+\left\| p_2\right\|_{\mathcal{H}_2}\right)\left\| \Psi\right\|_{\Z},
\end{eqnarray}
with $\displaystyle C_{\mathcal{G}_\Gamma}:=\max\left\{\frac{dC_{0,1}}{2\mu}, \; \frac{dC_{0,2}}{2\mu}\right\}.$

Finally, we estimate the third term in the right-hand side of equation~\eqref{eq:continuity}. Taking into account the Cauchy--Schwarz inequality, the trace theorem (Lemma~\ref{thm:tracetheorem}), the trace theorem on the tangential derivative (Lemma~\ref{thm:tracetheoremtangential}) and the Sobolev embedding on the interface (Lemma~\ref{lem:embedding-interface}), we get 
\begin{eqnarray}
\label{eq:continuity6}
      |\mathcal{F}_\Gamma(p_1,p_2; \vec{W})| &\le& \left\| p_1\right\|_{L^2(\Gamma)} \left\|W_{\vec{n}}\right\|_{L^2(\Gamma)} +\left\| p_2\right\|_{L^2(\Gamma)} \left\|W_{\vec{n}}\right\|_{L^2(\Gamma)} \nonumber
     \\
     &\quad&+  \frac{ \left(2\mu \ell^3 N_{2,1}^\mathrm{bl}+d\ell^2\right) M_{1,1}^{1,\mathrm{bl}}  }{\frac{1}{6}d^2 + \frac{2}{3}d \mu \ell  N_{1,1}^\mathrm{bl} +\frac{2}{3} d\mu \ell N_{2,1}^\mathrm{bl}  + 2 \mu^2 \ell^2  N_{1,1}^\mathrm{bl}  N_{2,1}^\mathrm{bl} }\left\| \frac{\partial p_1}{\partial \vec{\tau}}\right\|_{H^{-1/2}(\Gamma)}\left\|W_{\vec{\tau}}\right\|_{H^{1/2}(\Gamma)} \nonumber \\
        &\quad&+\frac{ \left(2\mu \ell^3  N_{1,1}^\mathrm{bl} +d\ell^2\right) M_{2,1}^{1,\mathrm{bl}} }{\frac{1}{6}d^2 + \frac{2}{3}d \mu \ell  N_{1,1}^\mathrm{bl} +\frac{2}{3} d\mu \ell N_{2,1}^\mathrm{bl}  + 2 \mu^2 \ell^2  N_{1,1}^\mathrm{bl}  N_{2,1}^\mathrm{bl} } \left\| \frac{\partial p_2}{\partial \vec{\tau}}\right\|_{H^{-1/2}(\Gamma)} \left\|W_{\vec{\tau}}\right\|_{H^{1/2}(\Gamma)}\nonumber\\
        &\le& C_{0,1}\left\| p_1\right\|_{\mathcal{H}_1} \left\|\vec{W}\right\|_{\mathcal{H}_\Gamma} +C_{0,2}\left\| p_2\right\|_{\mathcal{H}_2} \left\|W_{\vec{n}}\right\|_{\mathcal{H}_\Gamma}\nonumber
     \\
     &\quad&+  \frac{ \left(2\mu \ell^3 N_{2,1}^\mathrm{bl}+d\ell^2\right) M_{1,1}^{1,\mathrm{bl}}   C_{\vec{\tau},1} \Ce}{\frac{1}{6}d^2 + \frac{2}{3}d \mu \ell  N_{1,1}^\mathrm{bl} +\frac{2}{3} d\mu \ell N_{2,1}^\mathrm{bl}  + 2 \mu^2 \ell^2  N_{1,1}^\mathrm{bl}  N_{2,1}^\mathrm{bl} }\left\| p_1\right\|_{\mathcal{H}_1} \left\|\vec{W}\right\|_{\mathcal{H}_\Gamma}\nonumber \\
        &\quad&+\frac{ \left(2\mu \ell^3  N_{1,1}^\mathrm{bl} +d\ell^2\right) M_{2,1}^{1,\mathrm{bl}} C_{\vec{\tau},2} \Ce}{\frac{1}{6}d^2 + \frac{2}{3}d \mu \ell  N_{1,1}^\mathrm{bl} +\frac{2}{3} d\mu \ell N_{2,1}^\mathrm{bl}  + 2 \mu^2 \ell^2  N_{1,1}^\mathrm{bl}  N_{2,1}^\mathrm{bl} } \left\| p_2\right\|_{\mathcal{H}_2} \left\|\vec{W}\right\|_{\mathcal{H}_\Gamma}\nonumber\\
        &\le&C_{\mathcal{F}_\Gamma} \left(\left\| p_1\right\|_{\mathcal{H}_1} +\left\| p_2\right\|_{\mathcal{H}_2} \right)\left\|\vec{W}\right\|_{\mathcal{H}_\Gamma},
\end{eqnarray}
with 
\begin{eqnarray*}
    C_{\mathcal{F}_\Gamma} := \max\left\{  C_{0,1}+ \frac{ \left(2\mu \ell^3 N_{2,1}^\mathrm{bl}+d\ell^2\right) M_{1,1}^{1,\mathrm{bl}}   C_{\vec{\tau},1} \Ce}{\frac{1}{6}d^2 + \frac{2}{3}d \mu \ell  N_{1,1}^\mathrm{bl} +\frac{2}{3} d\mu \ell N_{2,1}^\mathrm{bl}  + 2 \mu^2 \ell^2  N_{1,1}^\mathrm{bl}  N_{2,1}^\mathrm{bl} },\right. \hspace{+30ex}\\
    \left.C_{0,2} +\frac{ \left(2\mu \ell^3  N_{1,1}^\mathrm{bl} +d\ell^2\right) M_{2,1}^{1,\mathrm{bl}} C_{\vec{\tau},2} \Ce}{\frac{1}{6}d^2 + \frac{2}{3}d \mu \ell  N_{1,1}^\mathrm{bl} +\frac{2}{3} d\mu \ell N_{2,1}^\mathrm{bl}  + 2 \mu^2 \ell^2  N_{1,1}^\mathrm{bl}  N_{2,1}^\mathrm{bl} } \right\}.
\end{eqnarray*}

Combining estimates~\eqref{eq:continuity1}--\eqref{eq:continuity6}, we show
continuity of the bilinear form $\mathcal{A}_\HYB$ in~\eqref{eq:weakformCoupledHybrid}:
\begin{eqnarray*}
    \left|\mathcal{A}_\HYB \left( \vec{\zeta} ; \vec{\theta}\right)\right| 
      &\le&  C_{\mathcal{A}_\HYB}\left(\| p_1\|_{\mathcal{H}_1} +\| p_2\|_{\mathcal{H}_2} +\| \vec{V}\|_{\mathcal{H}_\Gamma} +\|P\|_{\Z} \right)\left(\| \varphi_1\|_{\mathcal{H}_1} +\| \varphi_2\|_{\mathcal{H}_2} +\| \vec{W}\|_{\mathcal{H}_\Gamma} +\|\Psi\|_{\Z} \right)\\
     &\le& C_{\mathcal{A}_\HYB} \|\vec \zeta\|_{\mathcal{H}_\HYB} \|\vec \theta \|_{\mathcal{H}_\HYB},
\end{eqnarray*}
where $\displaystyle C_{\mathcal{A}_\HYB}:=\max \left\{ C_{\mathcal{A}_\PM}, \; C_{\mathcal{A}_\Gamma}, \; d, \; \frac{d}{\mu}, \; C_{\mathcal{F}_{\gamma_1,\gamma_2}},\;C_{\mathcal{G}_{\Gamma}},\;C_{\mathcal{F}_{\Gamma}}
\right\}$.

To finalise the proof, we check coercivity of the bilinear operator $\mathcal{A}_\HYB$ in~\eqref{eq:weakformCoupledHybrid}:
\begin{eqnarray}
    \mathcal{A}_\HYB (\vec{\zeta}; \vec{\zeta}  )&=&  \mathcal{A}_\Gamma\left( \vec{V}; \vec{V}\right) + \mathcal{E}_\Gamma\left(P;P \right)    +\mathcal{A}_\PM  (p_1,p_2; p_1, p_2)\nonumber  \\
      &\quad&-  \mathcal{F}_\Gamma( p_1, p_2;\vec{V})
        + \mathcal{G}_\Gamma (p_1, p_2; P )
     + \mathcal{A}_{\gamma_1,\gamma_2} (p_1,p_2; p_1, p_2)+\mathcal{F}_{\gamma_1, \gamma_2}(\vec{V}, P ; p_1, p_2) \nonumber\\
   &=& \underbrace{
   \sum_{i=1,2} \int_{\Omega_i} \left( \frac{\ten{K}_i}{\mu} \nabla p_i\right) \cdot \nabla p_i~\mathrm{d}\vec{x}}_{T_1}
  \nonumber\\
  &\quad& + \underbrace{d\mu \left\|\frac{\partial \vec{V}}{\partial \vec{\tau}}\right\|^2_{L^2(\Gamma)}  +\frac{\left(2 \mu d+2\mu^2 \ell N_{1,1}^\mathrm{bl}+2\mu^2 \ell N_{2,1}^\mathrm{bl}  \right)}{\frac{1}{6}d^2 + \frac{2}{3}d \mu \ell  N_{1,1}^\mathrm{bl} + \frac{2}{3}d\mu \ell N_{2,1}^\mathrm{bl}  + 2 \mu^2 \ell^2  N_{1,1}^\mathrm{bl}  N_{2,1}^\mathrm{bl} }  \left\| V_\vec{\tau}\right\|^2_{L^2(\Gamma)}}_{T_2} \nonumber \\
  &\quad&+\underbrace{\frac{ \left(2\mu \ell^3 N_{2,1}^\mathrm{bl}+d\ell^2\right)\MI}{\frac{1}{6}d^2 + \frac{2}{3}d \mu \ell  N_{1,1}^\mathrm{bl} +\frac{2}{3} d\mu \ell N_{2,1}^\mathrm{bl}  + 2 \mu^2 \ell^2  N_{1,1}^\mathrm{bl}  N_{2,1}^\mathrm{bl} }\int_\Gamma \left(\nabla p_1  \cdot \vec{\tau}|_{\gamma_1}\right) V_{\vec{\tau}}~\mathrm{d} s}_{T_3} \nonumber\\
 &\quad&+\underbrace{\frac{ \left(2\mu \ell^3  N_{1,1}^\mathrm{bl} +d\ell^2\right)\MII }{\frac{1}{6}d^2 + \frac{2}{3}d \mu \ell  N_{1,1}^\mathrm{bl} +\frac{2}{3} d\mu \ell N_{2,1}^\mathrm{bl}  + 2 \mu^2 \ell^2  N_{1,1}^\mathrm{bl}  N_{2,1}^\mathrm{bl} } \int_\Gamma \left(\nabla p_2  \cdot \vec{\tau}|_{\gamma_2}\right) V_{\vec{\tau}}~\mathrm{d} s}_{T_4}\nonumber \\
 &\quad& +\frac{d}{\mu} \|P\|^2_{L^2(\Gamma)}+\frac{d}{3\mu} \left(\|p_1|_{\gamma_1}\|_{L^2(\Gamma)}^2 + \|p_2|_{\gamma_2}\|_{L^2(\Gamma)}^2\right) \nonumber\\
 &\quad&-\frac{d}{\mu}\int_\Gamma  \left(p_1 |_{\gamma_1}+p_2 |_{\gamma_2}\right) P~\mathrm{d} s +  \frac{d}{3\mu} \int_{\Gamma}  (p_1|_{\gamma_1}) (p_2|_{\gamma_2}) ~\mathrm{d} s.
       \label{eq:CoecivityBiOperater}
\end{eqnarray}
For the last four terms in the right-hand side of~\eqref{eq:CoecivityBiOperater} it holds
\begin{eqnarray}
\label{eq:est-P}
   &\quad&\|P\|^2_{L^2(\Gamma)}+\frac{1}{3} \left(\|p_1|_{\gamma_1}\|_{L^2(\Gamma)}^2 + \|p_2|_{\gamma_2}\|_{L^2(\Gamma)}^2\right)-\int_\Gamma  \left(p_1 |_{\gamma_1}+p_2 |_{\gamma_2}\right) P~\mathrm{d} s +  \frac{1}{3} \int_{\Gamma}  (p_1|_{\gamma_1}) (p_2|_{\gamma_2}) ~\mathrm{d} s \nonumber\\
   &=& \left\|P - \frac{1}{2}\left(p_1|_{\gamma_1} + p_2 |_{\gamma_2}\right)\right\|_{L^2(\Gamma)}^2 + \frac{1}{12}\|p_1|_{\gamma_1} - p_2 |_{\gamma_2}\|_{L^2(\Gamma)}^2.
\end{eqnarray}
Using the fact, that the permeability tensor is symmetric, positive definite and bounded~\eqref{equ:boundedK} and applying the Poincar\'{e} theorem (Lemma~\ref{thm:poicare}), we obtain the following estimate for the term $T_1$ in~\eqref{eq:CoecivityBiOperater}:
\begin{align}
\label{eq:est-Apm}
    \sum_{i=1,2} \int_{\Omega_i} \left( \frac{\ten{K}_i}{\mu} \nabla p_i\right) \cdot \nabla p_i~\mathrm{d}\vec{x} \geq \sum_{i\in \{1,2\}} \frac{k_{\min,i}}{\mu \tilde{C}_{P,i}} \|p_i\|^2_{\mathcal{H}_i}.
\end{align}
For the term $T_2$, we get
\begin{eqnarray}
\label{eq:est-Ev-orig}
    d\mu \left\|\frac{\partial \vec{V}}{\partial \vec{\tau}}\right\|^2_{L^2(\Gamma)}  +B(d)  \left\| V_\vec{\tau}\right\|^2_{L^2(\Gamma)} \geq d\mu \left\|\frac{\partial V_{\vec{\tau}}}{\partial \vec{\tau}}\right\|^2_{L^2(\Gamma)}  +B(d) \left\| V_\vec{\tau}\right\|^2_{L^2(\Gamma)},
    % &\geq& \frac{1}{C_I}\min\bigg\{\frac{\left(2 \mu d+2\mu^2 \ell N_{1,1}^\mathrm{bl}+2\mu^2 \ell N_{2,1}^\mathrm{bl}  \right)}{\frac{1}{6}d^2 + \frac{2}{3}d \mu \ell  N_{1,1}^\mathrm{bl} + \frac{2}{3}d\mu \ell N_{2,1}^\mathrm{bl}  + 2 \mu^2 \ell^2  N_{1,1}^\mathrm{bl}  N_{2,1}^\mathrm{bl} }, \\
    % &\quad& \hspace{9ex}\sqrt{2d\mu \frac{\left(2 \mu d+2\mu^2 \ell N_{1,1}^\mathrm{bl}+2\mu^2 \ell N_{2,1}^\mathrm{bl}  \right)}{\frac{1}{6}d^2 + \frac{2}{3}d \mu \ell  N_{1,1}^\mathrm{bl} + \frac{2}{3}d\mu \ell N_{2,1}^\mathrm{bl}  + 2 \mu^2 \ell^2  N_{1,1}^\mathrm{bl}  N_{2,1}^\mathrm{bl} } }\bigg\} \|\vec{V}_{\vec{\tau}}\|_{H^{1/2}(\Gamma)}^2,
\end{eqnarray}
with $B(d)$ defined in~\eqref{eq:abbr-B}.
We can estimate it further, using the Sobolev interpolation inequality (Lemma~\ref{lem:Sobolev-interpolation}):
\begin{eqnarray}
\label{eq:est-Ev}
    &\quad&\left(d\mu \left\|\frac{\partial V_{\vec{\tau}}}{\partial \vec{\tau}}\right\|^2_{L^2(\Gamma)}  + B(d)  \left\| V_\vec{\tau}\right\|^2_{L^2(\Gamma)}\right)^2 \nonumber\\
    &\geq& B(d)^2\left\| V_\vec{\tau}\right\|^4_{L^2(\Gamma)}  + 2d\mu B(d)\left\|\frac{\partial V_{\vec{\tau}}}{\partial\vec{\tau}}\right\|^2_{L^2(\Gamma)} \left\|V_{\vec{\tau}}\right\|_{L^2(\Gamma)}^2 \nonumber\\
    &\geq&  \min\left\{B(d)^2, \,2d\mu B(d) \right\}\|V_{\vec{\tau}}\|_{L^2(\Gamma)}^2\left(\|V_{\vec{\tau}}\|_{L^2(\Gamma)}^2 + \left\|\frac{\partial V_{\vec{\tau}}}{\partial \vec{\tau}}\right\|^2_{L^2(\Gamma)}\right) \nonumber\\
    &\geq& C^2_{T_2}(d) \|V_{\vec{\tau}}\|_{H^{1/2}(\Gamma)}^4, 
    % &\geq&  \frac{1}{C_I^2}\min\{B(d)^2, \, 2d\mu B(d)\}\|\vec{V}_{\vec{\tau}}\|_{H^{1/2}(\Gamma)}^4.
\end{eqnarray}
with $B(d)$ and $C_{T_2}(d)$ defined in~\eqref{eq:abbr}.
Taking into account the trace theorem on the tangential derivative (Lemma~\ref{thm:tracetheoremtangential}) and the generalised Young inequality $ 2ab \le a^2/\delta + \delta b^2\textnormal{ for } \delta>0,$  we obtain for the terms $T_3$ and $T_4$ in~\eqref{eq:CoecivityBiOperater} with $b_i$ given in~\eqref{eq:abbr}:
\begin{eqnarray*}
    |T_3 + T_4| &\leq& b_1(d)\left\|\frac{\partial p_1}{\partial \vec{\tau}} \right\|_{H^{-1/2}(\Gamma)} \|V_{\vec{\tau}}\|_{H^{1/2}(\Gamma)} +b_2(d)\left\|\frac{\partial p_2}{\partial \vec{\tau}} \right\|_{H^{-1/2}(\Gamma)} \|V_{\vec{\tau}}\|_{H^{1/2}(\Gamma)} \nonumber\\
    &\leq& b_1(d) C_{\vec{\tau},1}\left\| p_1 \right\|_{\mathcal{H}_1} \|V_{\vec{\tau}}\|_{H^{1/2}(\Gamma)} +b_2(d) C_{\vec{\tau},2}\left\| p_2 \right\|_{\mathcal{H}_2} \|V_{\vec{\tau}}\|_{H^{1/2}(\Gamma)} \nonumber
    \\ &\leq& \frac{b_1(d)}{2}\left(\frac{C_{\vec{\tau},1}^2}{\delta_1}\left\| p_1 \right\|^2_{\mathcal{H}_1} +\delta_1\|V_{\vec{\tau}}\|_{H^{1/2}(\Gamma)}^2\right)+\frac{b_2(d)}{2}\left(\frac{C_{\vec{\tau},2}^2}{\delta_2}\left\| p_2 \right\|^2_{\mathcal{H}_2} +\delta_2\|V_{\vec{\tau}}\|_{H^{1/2}(\Gamma)}^2\right),
\end{eqnarray*}
where the following parameters are chosen
\begin{align*}
    \delta_i = \frac{b_i(d) \mu C_{\vec{\tau},i}^2 \tilde{C}_{P,i}}{2\theta k_{\min,i}}, \quad \theta \in (0,1), \qquad i \in \{1,2\}.
\end{align*}
This leads to
\begin{eqnarray*}
    |T_3 + T_4| &\leq& \theta \sum_{i \in \{1,2\}} \frac{k_{\min,i}}{\mu \tilde{C}_{P,i}} \|p_i\|_{\mathcal{H}_i}^2 + \frac{1}{4\theta}\left(\sum_{i \in \{1,2\}} \frac{b_i(d)^2\mu C_{\vec{\tau},i}^2 \tilde{C}_{P,i}}{k_{\min,i}}\right) \|V_{\vec{\tau}}\|^2_{H^{1/2}(\Gamma)}.
\end{eqnarray*}
Taking into account the estimates in~\eqref{eq:est-Ev} and~\eqref{eq:est-Ev-orig}, we obtain
\begin{eqnarray}
\label{eq:est-T3-T4}
    |T_3 + T_4| &\leq& \theta \sum_{i \in \{1,2\}} \frac{k_{\min,i}}{\mu \tilde{C}_{P,i}} \|p_i\|_{\mathcal{H}_i}^2  + \frac{1}{4\theta}\frac{C_{T_{3,4}}(d)}{C_{T_2}(d)} \left(d\mu \left\|\frac{\partial \vec{V}}{\partial \vec{\tau}}\right\|^2_{L^2(\Gamma)}  + B(d)  \left\| V_\vec{\tau}\right\|^2_{L^2(\Gamma)}\right),
\end{eqnarray}
where $C_{T_{3,4}}(d)$ is defined in~\eqref{eq:abbr-CT34}.
We assume
\begin{align}
\label{eq:assumption}
    \frac{C_{T_{3,4}}(d)}{4C_{T_2}(d)} < 1
\end{align}
and choose
\begin{align}
\label{eq:theta}
   \theta =  \sqrt{\frac{C_{T_{3,4}}(d)}{4C_{T_2}(d)}}.
\end{align}
\begin{remark}
Note that all physical parameters in~\eqref{eq:assumption} are positive. 
Moreover, we analysed~\eqref{eq:assumption} for different physical parameters and concluded that this assumption is satisfied for a wide range of applications.
\end{remark}

We estimate~\eqref{eq:CoecivityBiOperater} by applying~\eqref{eq:est-P},~\eqref{eq:est-Apm} and~\eqref{eq:est-T3-T4} together with the Poincar\'{e} inequality (Lemma~\ref{thm:poicare}):
\begin{eqnarray}
    \mathcal{A}_\HYB (\vec{\zeta}; \vec{\zeta}) &\geq&(1-\theta) \left(\left(d\mu \left\|\frac{\partial \vec{V}}{\partial \vec{\tau}}\right\|^2_{L^2(\Gamma)}  + B(d)  \left\| V_\vec{\tau}\right\|^2_{L^2(\Gamma)}\right) + \sum_{i \in \{1,2\}} \frac{k_{\min,i}}{\mu \tilde{C}_{P,i}} \|p_i\|^2_{\mathcal{H}_i}\right) \nonumber\\
    &\quad& + \frac{d}{\mu}\left\|P - \frac{1}{2}\left(p_1|_{\gamma_1} + p_2 |_{\gamma_2}\right)\right\|_{L^2(\Gamma)}^2 + \frac{d}{12\mu}\|p_1|_{\gamma_1} - p_2 |_{\gamma_2}\|_{L^2(\Gamma)}^2 \nonumber\\
    &\geq& \frac{(1-\theta) d\mu }{\tilde{C}_{P,\Gamma}} \|\vec{V}\|^2_{\mathcal{H}_{\Gamma}} + \sum_{i \in \{1,2\}} (1-\theta)\frac{k_{\min,i}}{\mu \tilde{C}_{P,i}} \|p_i\|^2_{\mathcal{H}_i}  + \frac{d}{\mu}\left\|P - \frac{1}{2}\left(p_1|_{\gamma_1} + p_2 |_{\gamma_2}\right)\right\|_{L^2(\Gamma)}^2.
\end{eqnarray}
Using the estimation
\begin{eqnarray*}
    \|P\|_{L^2(\Gamma)}^2 &\leq& 2 \left\|P - \frac{1}{2}\left(p_1|_{\gamma_1} + p_2 |_{\gamma_2}\right)\right\|_{L^2(\Gamma)}^2 + \|p_1|_{\gamma_1}\|_{L^2(\Gamma)}^2 +\|p_2|_{\gamma_2}\|_{L^2(\Gamma)}^2,
\end{eqnarray*}
and applying the trace theorem (Lemma~\ref{thm:tracetheorem}), we obtain
\begin{eqnarray}
    \label{eq:est-P-2}
        \|P\|_{L^2(\Gamma)}^2 &\leq& 2 \left\|P - \frac{1}{2}\left(p_1|_{\gamma_1} + p_2 |_{\gamma_2}\right)\right\|_{L^2(\Gamma)}^2 + C_{0,1}^2\|p_1\|^2_{\mathcal{H}_1} + C_{0,2}^2\|p_2\|^2_{\mathcal{H}_2}.
\end{eqnarray}
Taking~\eqref{eq:est-P-2} into account, we conclude
\begin{align}
    \mathcal{A}_\HYB (\vec{\zeta}; \vec{\zeta}) \geq C^*_{\mathcal{A}_\HYB} \| \vec{\zeta}\|^2_{\mathcal{H}_\HYB},
\end{align}
with 
\begin{align*}
    C^*_{\mathcal{A}_\HYB} := \min \left\{\frac{(1-\theta) d\mu }{\tilde{C}_{P,\Gamma}},\, \frac{(1-\theta)k_{\min,1}}{(1+C_{0,1}^2)\mu \tilde{C}_{P,1}},\, \frac{(1-\theta)k_{\min,2}}{(1+C_{0,2}^2)\mu \tilde{C}_{P,2}},\, \frac{d}{2\mu}\right\}.
\end{align*}
This completes the proof.

\section{Numerical Simulations}\label{sec:NumSimulation}
In this section, we study the developed hybrid-dimensional fracture model numerically. The model is discretised using the multipoint flux approximation (MPFA) scheme in the porous-medium domains $\Omega_1$ and $\Omega_2$ and the second-order finite difference method on the interface $\Gamma$. The problem is then solved monolithically using our C++ code. First, we confirm numerically the second order of convergence for the discrete hybrid-dimensional model in Section~\ref{sec:convergence-order}. Then, in Section~\ref{sec:validation} we validate the developed model against the pore-scale resolved model and compare it with two other hybrid-dimensional models developed in~\cite{Rybak-Metzger-20}.

\subsection{Convergence order}
\label{sec:convergence-order}
Here, we demonstrate the second-order convergence considering the following test case~\ref{test:1}.
\begin{test}
\label{test:1}
    The coupled flow domain $\Omega$ comprises the porous-medium regions $\Omega_1 = [0,1.0] \times [0.0+0.5d, 1.0+0.5d],\, \Omega_2~=~[0,1.0] \times [-1.0-0.5d, 0.0-0.5d]$ and the lower-dimensional complex interface $\Gamma = [0,1.0] \times \{0.0\}$ (fracture). Both porous media are isotropic with permeability tensors $\ten K_i = k_i \ten I, \, k_i > 0$ for $i=1,2$. The model parameters are $\mu = 1.0$ and $k_i = 10^{-2}, \, i=1,2$. For this test case, the boundary layer constants are set to $N_{i,1}^{\mathrm{bl}} = M_{i,1}^{1,\mathrm{bl}} = 1.0,\, i=1,2$. The analytical solution for the proposed hybrid-dimensional model is given by
    \begin{align*}
        p_1 &= -100(y-y_{\gamma_1})\sin(x)\exp(d), \quad p_2 = -50(y-y_{\gamma_2})^2\sin(x)\exp(d), \\
        V_{\vec{\tau}} &= \frac{1}{d}\exp(d)\cos(x), \quad V_\vec{n} = \frac{1}{2}\exp(d)\sin(x), \quad P = \frac{1}{d} \exp(d) \sin(x),
    \end{align*}
    where $y_{\gamma_i}$ refers to the $y$-value of the interfaces $\gamma_i, \; i=1,2$.

    We choose Dirichlet boundary conditions on the outer boundary of the flow domain $\partial \Omega$. The right-hand sides are obtained by substituting the analytical solution and the physical parameters into~\eqref{eq:Darcy-law},~\eqref{equ:AveragedMass},~\eqref{eqn:averagedMomentNormalFinal},~\eqref{eqn:averagedMomentTangentialFinal}, \eqref{eq:transmission1} and~\eqref{eq:transmission2}.
\end{test}
For the convergence analysis, we compute the relative $L_2$-errors between the numerical simulation results $f_h$ for grid width $h$ and the analytical solution $f$ for all primary variables
$$\varepsilon_f = \|f -f_h\|_{L_2(\Omega)}/ \|f\|_{L_2(\Omega)},\quad f\in \{p_1,\,p_2\} \quad \text{and} \quad \varepsilon_f = \|f -f_h\|_{L_2(\Gamma)},\quad f\in \{V_\vec{\tau}, V_\vec{n}, P\}.$$ Using these $L_2$-errors, we obtain the convergence rates of the discretisation scheme. We consider two different fracture apertures $d= 10^{-1}$ (Tab.~\ref{tab:L2-error-test-d=1e-1} and Tab.~\ref{tab:convergence-rates-test-d=1e-1}) and $d= 10^{-2}$ (Tab.~\ref{tab:L2-error-test-d=1e-2} and Tab.~\ref{tab:convergence-rates-test-d=1e-2}).
The results presented in Tab.~\ref{tab:L2-error-test-d=1e-1}--Tab.~\ref{tab:convergence-rates-test-d=1e-2} confirm the second-order convergence of the discretisation scheme. 

\begin{table}[!ht]
    \centering
    \begin{tabular}{|c||c|c|c|c|c|c|c|}
        \hline
        Grid 
        & $h=1/8$ 
        & $h=1/16$ 
        & $h=1/32$ 
        & $h=1/64$ 
        & $h=1/128$ 
        & $h=1/256$ 
        & $h=1/512$\\
        \hline \hline
        
        $\varepsilon_{p_1}$ 
        & $2.61 \times 10^{-3}$ & $6.52 \times 10^{-4}$ & $1.63 \times 10^{-4}$ & $4.06 \times 10^{-5}$ & $1.01 \times 10^{-5}$ & $2.54 \times 10^{-6}$ & $6.34 \times 10^{-7}$ \\
        \hline
        
        $\varepsilon_{p_2}$ 
        & $1.12 \times 10^{-2}$ & $3.18 \times 10^{-3}$ & $8.4 \times 10^{-4}$ & $2.15 \times 10^{-4}$ & $5.45 \times 10^{-5}$ & $1.37 \times 10^{-5}$ & $3.44 \times 10^{-6}$ \\
        \hline
        
        $\varepsilon_{V_{\vec{\tau}}}$ 
        & $3.43 \times 10^{-3}$ & $1.11 \times 10^{-3}$ & $3.24 \times 10^{-4}$ & $8.89 \times 10^{-5}$ & $2.33 \times 10^{-5}$ & $5.94 \times 10^{-6}$ & $1.5 \times 10^{-6}$ \\
        \hline
        
        $\varepsilon_{V_{\vec{n}}}$ 
        & $2.07 \times 10^{-2}$ & $6.11 \times 10^{-3}$ & $1.64 \times 10^{-3}$ & $4.23 \times 10^{-4}$ & $1.07 \times 10^{-4}$ & $2.69 \times 10^{-5}$ & $6.75 \times 10^{-6}$ \\
        \hline
        
        $\varepsilon_{P}$ 
        & $3.81 \times 10^{-2}$ & $1.0 \times 10^{-2}$ & $2.68 \times 10^{-3}$ & $7.15 \times 10^{-4}$ & $1.87 \times 10^{-4}$ & $4.8 \times 10^{-5}$ & $1.21 \times 10^{-5}$\\
        \hline
    \end{tabular}
    \caption{Relative errors for all primary variables in Test~\ref{test:1} for $d=10^{-1}$.}
    \label{tab:L2-error-test-d=1e-1}
\end{table}

\begin{table}[!ht]
    \centering
    \begin{tabular}{|c||c|c|c|c|c|c|}
        \hline
        Grid
        & $8/16$
        & $16/32$
        & $32/64$
        & $64/128$
        & $128/256$
        & $256/512$\\
        \hline \hline

        ${p_1}$
        & 2.00 & 2.00 & 2.01 & 2.01 & 1.99 & 2.00\\
        \hline

        ${p_2}$
        & 1.82 & 1.92 & 1.97 & 1.98 & 1.99 & 1.99 \\
        \hline

        ${V_{\vec{\tau}}}$
        & 1.63 & 1.78 & 1.87 & 1.93 & 1.97 & 1.99 \\
        \hline

        ${V_{\vec{n}}}$
        & 1.76 & 1.90 & 1.96 & 1.98 & 1.99 & 2.00 \\
        \hline

        ${P}$
        & 1.93 & 1.90 & 1.91 & 1.94 & 1.96 & 1.99 \\
        \hline
    \end{tabular}
    \caption{Convergence orders for all primary variables in Test~\ref{test:1} for $d=10^{-1}$.}
    \label{tab:convergence-rates-test-d=1e-1}
\end{table}

\begin{table}[!ht]
    \centering
    \begin{tabular}{|c||c|c|c|c|c|c|c|c|}
        \hline
        Grid 
        & $h=1/8$ 
        & $h=1/16$ 
        & $h=1/32$ 
        & $h=1/64$ 
        & $h=1/128$ 
        & $h=1/256$
        & $h=1/512$\\
        \hline \hline
        
        $\varepsilon_{p_1}$ 
        & $2.72 \times 10^{-3}$ & $6.81 \times 10^{-4}$ & $1.7 \times 10^{-4}$ & $4.26 \times 10^{-5}$ & $1.07 \times 10^{-5}$ & $2.67 \times 10^{-6}$ & $6.66 \times 10^{-7}$\\
        \hline
        
        $\varepsilon_{p_2}$ 
        & $1.12 \times 10^{-2}$ & $3.16 \times 10^{-3}$ & $8.34 \times 10^{-4}$ & $2.14 \times 10^{-4}$ & $5.41 \times 10^{-5}$ & $1.36 \times 10^{-5}$ & $3.41 \times 10^{-6}$\\
        \hline
        
        $\varepsilon_{V_{\vec{\tau}}}$ 
        & $1.89 \times 10^{-3}$ & $5.41 \times 10^{-4}$ & $1.46 \times 10^{-4}$ & $3.81 \times 10^{-5}$ & $9.72 \times 10^{-6}$ & $2.45 \times 10^{-6}$ & $6.15 \times 10^{-7}$\\
        \hline
        
        $\varepsilon_{V_{\vec{n}}}$ 
        & $2.55 \times 10^{-2}$ & $8.95 \times 10^{-3}$ & $2.55 \times 10^{-3}$ & $6.71 \times 10^{-4}$ & $1.71 \times 10^{-4}$ & $4.31 \times 10^{-5}$ & $1.08 \times 10^{-5}$\\
        \hline
        
        $\varepsilon_{P}$ 
        & $9.59 \times 10^{-3}$ & $2.9 \times 10^{-3}$ & $8.47 \times 10^{-4}$ & $2.37 \times 10^{-4}$ & $6.37 \times 10^{-5}$ & $1.65 \times 10^{-5}$ & $4.2 \times 10^{-6}$\\
        \hline
    \end{tabular}
    \caption{Relative errors for all primary variables in Test~\ref{test:1} for $d=10^{-2}$.}
    \label{tab:L2-error-test-d=1e-2}
\end{table}

\begin{table}[!ht]
    \centering
    \begin{tabular}{|c||c|c|c|c|c|c|}
        \hline
        Grid 1/Grid 2
        & $8/16$
        & $16/32$
        & $32/64$
        & $64/128$
        & $128/256$
        & $256/512$\\
        \hline \hline

        ${p_1}$
        & 2.00 & 2.00 & 2.00 & 1.99 & 2.00 & 2.00\\
        \hline

        ${p_2}$
        & 1.83 & 1.92 & 1.96 & 1.98 & 1.99 & 2.00\\
        \hline

        ${V_{\vec{\tau}}}$
        & 1.80 & 1.89 & 1.94 & 1.97 & 1.99 & 1.99\\
        \hline

        ${V_{\vec{n}}}$
        & 1.51 & 1.81 & 1.93 & 1.97 & 1.99 & 2.00\\
        \hline

        ${P}$
        & 1.73 & 1.78 & 1.84 & 1.90 & 1.95 & 1.97\\
        \hline
    \end{tabular}
    \caption{Convergence orders for all primary variables in Test~\ref{test:1} for $d=10^{-2}$.}
    \label{tab:convergence-rates-test-d=1e-2}
\end{table}

\subsection{Model comparison and validation}
\label{sec:validation}
\begin{figure}[!ht]
    \centering
    \includegraphics[width=0.7\linewidth]{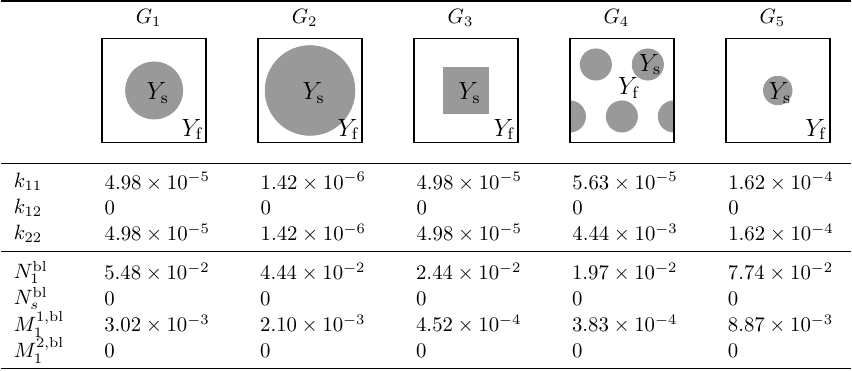}
    \caption{Permeability values and boundary layer constants for different pore geometries.}
    \label{fig:permeability-table}
\end{figure}

In this section, we validate the newly developed hybrid-dimensional Stokes--Darcy model for fractured porous media using pore-scale simulation results (DNS). We also compare the new model to previously developed hybrid-dimensional models, where instead of the generalised interface conditions (ER) the classical set of coupling conditions with either the Beavers--Joseph (BJ) or the Beavers--Joseph--Saffman (BJS) condition on the tangential velocity was taken into account~\cite{Rybak-Metzger-20}.
We validate the models by averaging the pore-scale velocity in the vertical direction across the fracture and comparing the resulting velocity profiles with those predicted by the hybrid-dimensional models.
% Due to oscillations in the averaged vertical pore-scale velocity, determining which model fits best here is infeasible and we limit this analysis to the horizontal velocity.
We consider three different flow scenarios (Section~\ref{sec:A},~\ref{sec:B} and~\ref{sec:C})
to show the applicability of the new hybrid-dimensional model to arbitrary flow directions in the vicinity of the interface. To investigate the influence of the pore geometry, we consider both circular and square inclusions as well as in-line and staggered arrangements, thereby covering different surface roughness characteristics in the vicinity of the fracture. The computed permeability and boundary layer values are presented in Fig.~\ref{fig:permeability-table}. The fluid viscosity is set to $\mu = 1.0$ and for the classical set of coupling conditions the Beavers--Joseph slip coefficient is chosen $\alpha_\BJ = 1.0$, following common practice in the literature. 

The pore-scale model~\eqref{eq:Stokes-micro} and~\eqref{eq:BC-micro}, the Stokes problem in the unit cell, e.g.~\cite{Hornung_97, Jaeger_Mikelic_96, Eggenweiler_Rybak_MMS20}, needed to derive the permeability values, and the boundary layer problems are solved using FreeFEM++ with Taylor--Hood finite elements~\cite{Hecht_12}. All hybrid-dimensional models are discretised with grid size $h= 2^{-9}$.

% \Ptodo{Surface roughness should be included in $\alpha_\BJ$. There is commonly nothing given on how to calculate the constant. Exception is the paper of Jäger \& Mikelic, which is a special case of the introduced generalised conditions. They can be interpreted as a modification of Beavers--Joseph conditions (see "A modification..."). }

\subsubsection{Flow scenario A}
\label{sec:A}
\begin{figure}[!ht]
    \centering
    \includegraphics[width=0.49\linewidth]{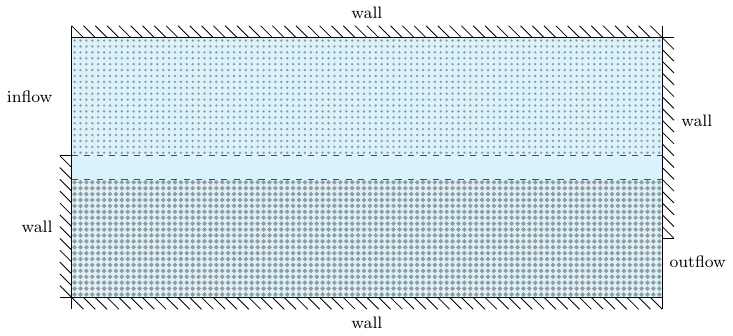}
    \raisebox{2.5mm}{
    \includegraphics[width=0.49\linewidth]{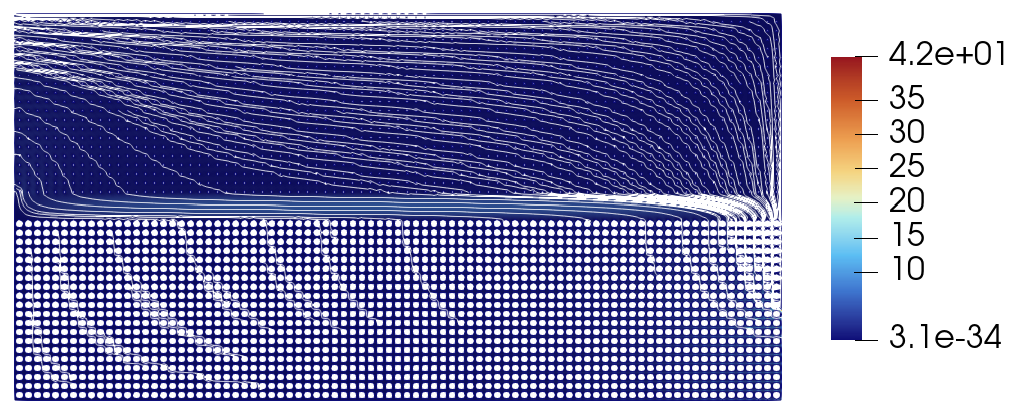}
    }
    \caption{Geometric setting and pore-scale velocity magnitude in \textit{Flow scenario A}.}
    \label{fig:geometric-setting-Test2}
\end{figure}

In the first test case, inflow is considered on the left boundary of the upper porous-medium domain $\Gamma_1^{\text{in}}$ and outflow on the right boundary of the bottom porous-medium domain $\Gamma_2^{\text{out}}$ (Fig.~\ref{fig:geometric-setting-Test2}).

% \begin{test}
% \label{test:2}
 We consider the flow domain $\Omega = [0, 5.0] \times [-1.1, 1.1]$, the inflow boundary $\Gamma_1^{\text{in}} = \{0\} \times [0.1, 5.0]$, the outflow boundary $\Gamma_2^{\text{out}} = \{5.0\} \times [-1.1, -0.6]$ and the remaining boundaries $\Gamma^{\text{wall}} = ([0, 5.0] \times \{1.1\}) \cup ([0, 5.0] \times \{-1.1\}) \cup (\{0\} \times [-1.1, 0.1]) \cup (\{5.0\} \times [-0.6, 1.1])$. In the hybrid-dimensional model, we have the porous-medium domains $\Omega_1 = [0,5.0] \times [0.1, 1.1]$ and $\Omega_2 = [0, 5.0] \times [-1.1, -0.1]$ and the complex interface $\Gamma = [0, 5.0] \times \{0\}$. We use the following boundary conditions
\begin{align}
    \vec{v} = (1,0) \text{ on } \Gamma_1^{\text{in}}, \quad p = 0 \text{ on } \Gamma_2^{\text{out}}, \quad \vec{v} = \textbf{0} \text{ on } \Gamma^{\text{wall}}.
\end{align}
% \end{test}

\begin{figure}[!ht]
    \centering
    \includegraphics[width=0.55\linewidth]{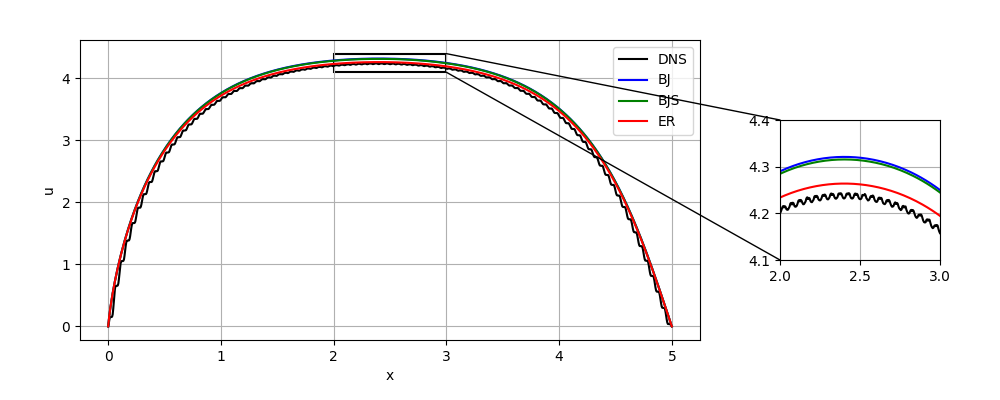} 
     \includegraphics[width=0.4\linewidth]{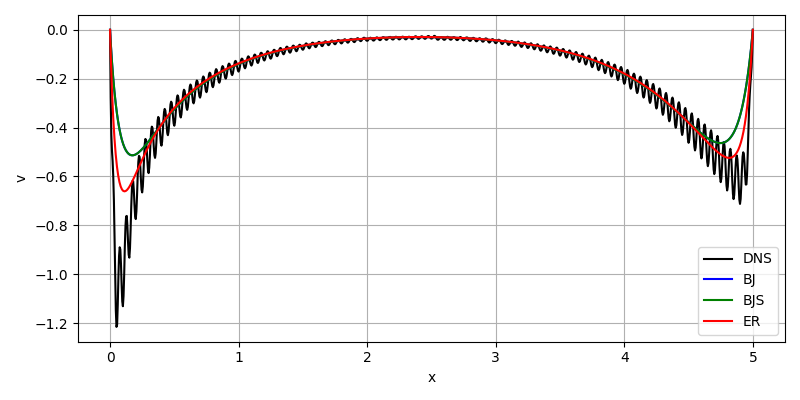} 
    \caption{Tangential (left) and normal (right) velocity component along the fracture in \textit{Flow scenario~A} with geometry $G_1$ in both porous-medium domains.}
    \label{fig:Test2-G1}
\end{figure}

\begin{figure}[!ht]
    \centering
   \includegraphics[width=0.55\linewidth]{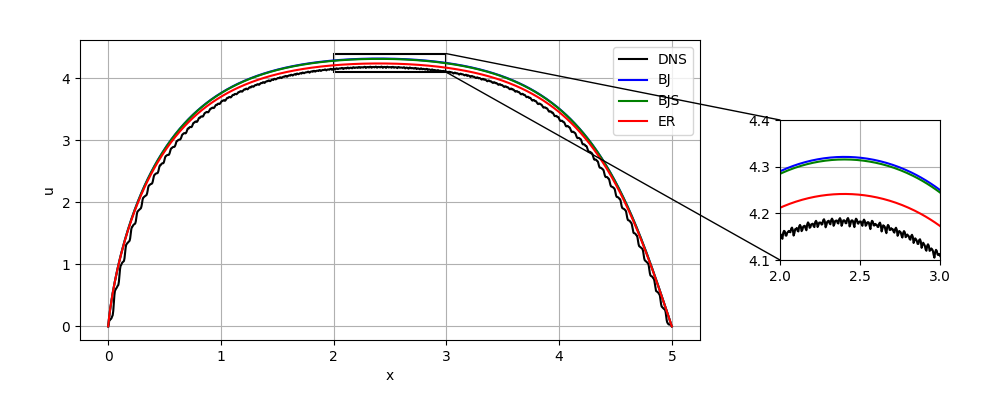}
     \includegraphics[width=0.4\linewidth]{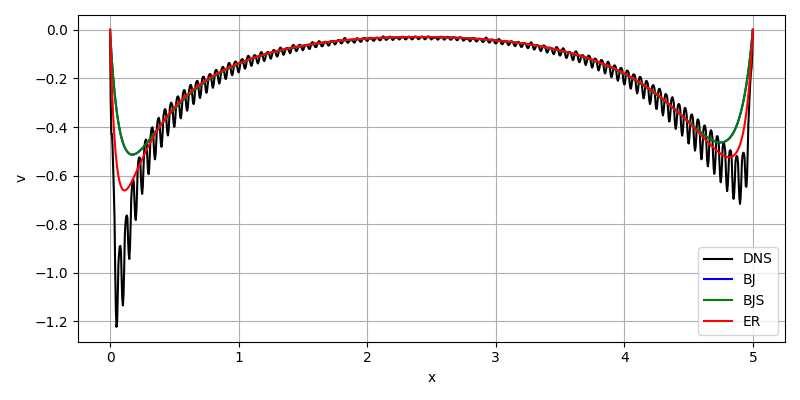} 
    \caption{Tangential (left) and normal (right) velocity component along the fracture in \textit{Flow scenario~A} with geometry $G_3$ in both porous-medium domains.}
    \label{fig:Test2-G3}
\end{figure}

% \begin{test}
% \label{test:2}
%  We consider the flow domain $\Omega = [0, 29.75] \times [-7.5, 7.5]$, the inflow boundary \ps{$\Gamma_1^{\text{in}} = \{0\} \times [0.5, 7.5]$}, the outflow boundary $\ps{\Gamma_2^{\text{out}} = \{29.75\} \times [-7.5, -4]}$ and the remaining boundaries $\Gamma^{\text{wall}} = ([0, 29.75] \times \{7.5\}) \cup ([0, 29.75] \times \{-7.5\}) \cup (\{0\} \times [-7.5, 0.5]) \cup (\{29.75\} \times [-4, 7.5])$. In the hybrid-dimensional model, we have the porous-medium domains $\Omega_1 = [0,29.75] \times [0.5, 7.5]$ and $\Omega_2 = [0, 29.75] \times [-7.5, -0.5]$ and the complex interface $\Gamma = [0, 29.75] \times \{0\}$. We use the following boundary conditions
% \begin{align}
%     \vec{v} = (1,0) \text{ on } \Gamma_1^{\text{in}}, \quad p = 0 \text{ on } \Gamma_2^{\text{out}}, \quad \vec{v} = \textbf{0} \text{ on } \Gamma^{\text{wall}}.
% \end{align}
% \end{test}
Below, we study the influence of the pore structure on the flow profile. First, we consider two different isotropic pore geometries with the same permeability, but different surface roughness characteristics which leads to different boundary layer constants (Fig.~\ref{fig:permeability-table}, geometries $G_1$ and $G_3$). In the first case, both porous-medium domains are constructed by $100 \times 20$ periodically distributed circular inclusions ($\ell = 0.05$) with radius $r=0.25 \ell$. In the second case, the porous media are build out of $100 \times 20$ square inclusions with length $a=0.2154\ell$.
The permeability values and boundary layer coefficients are provided in Fig.~\ref{fig:permeability-table}. The corresponding velocity plots are displayed in Fig.~\ref{fig:Test2-G1} and Fig.~\ref{fig:Test2-G3}. For both geometries, the newly developed hybrid-dimensional model (ER) is closer to the pore-scale resolved model (DNS) than the previously developed ones (BJ and BJS). The discrepancy between the two approaches is more pronounced for the square inclusions (geometry $G_3$).

\begin{figure}[!ht]
    \centering
    \includegraphics[width=0.55\linewidth]{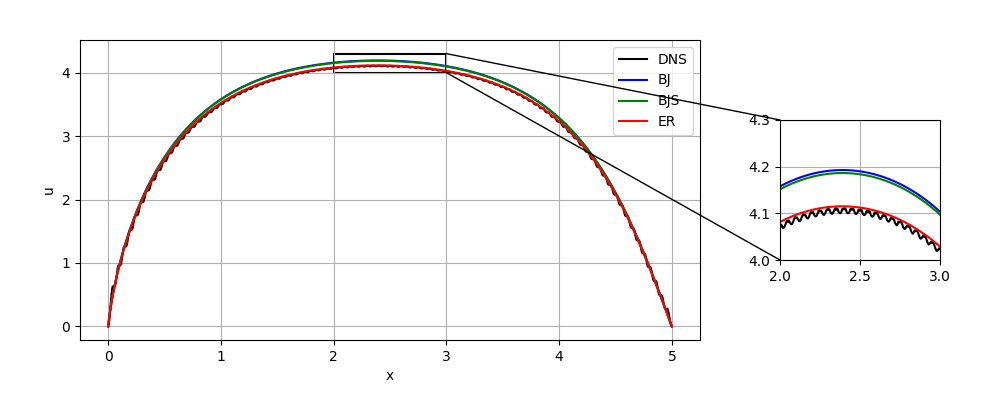}
    \includegraphics[width=0.4\linewidth]{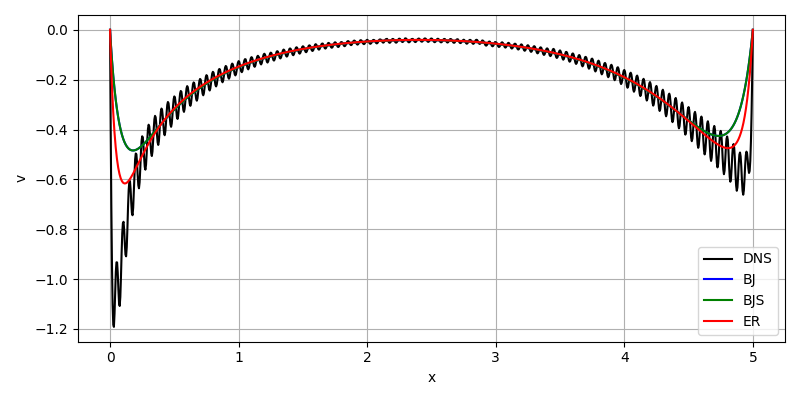}
    \caption{Tangential (left) and normal  (right) velocity component along the fracture in \textit{Flow scenario~A} with geometry $G_4$ in both porous-medium domains.}
    \label{fig:Test2-G4}
\end{figure}

Next, we expand our investigation to orthotropic porous media. In both porous media, we consider $100 \times 20$ circular solid inclusions with the same surface roughness as geometry $G_1$ and arrange them in a staggered manner (Fig.~\ref{fig:permeability-table}, geometry $G_4$). This configuration corresponds to a characteristic pore size of $\ell = 0.1$ and the radius $r = 0.125\ell$ for the circular inclusions. The velocity profiles shown in Fig.~\ref{fig:Test2-G4} once again confirm the superior agreement of the newly developed hybrid-dimensional model with the pore-scale simulations.

Finally, we consider a case where we have isotropic porous media with vastly different permeability in the two porous-medium domains (Fig.~\ref{fig:permeability-table}, geometries $G_2$ and $G_5$). The top region $\Omega_1$ is constructed by $100 \times 20$ circular inclusions ($\ell = 0.05$) with radius $r = 0.125 \ell$, which leads to a high permeability $k_{11} = k_{22} = 1.62\cdot 10^{-4}$. In the bottom porous-medium domain $\Omega_2$, we also have $100 \times 20$ circular inclusions constructed with a larger inclusion radius, $r= 0.437 \ell$, which yields low permeability ($k_{11} = k_{22} = 1.42 \cdot 10^{-6}$). The new hybrid-dimensional model shows again good agreement with the pore-scale simulations, while the previously developed models (BJ and BJS) from~\cite{Rybak-Metzger-20} are not as accurate. 

To conclude, for \textit{Flow scenario A}, the newly developed hybrid-dimensional model consistently agrees more closely with the pore-scale simulations than the models developed in~\cite{Rybak-Metzger-20} across all considered pore geometries. 

\begin{figure}[!ht]
    \centering
    \includegraphics[width=0.55\linewidth]{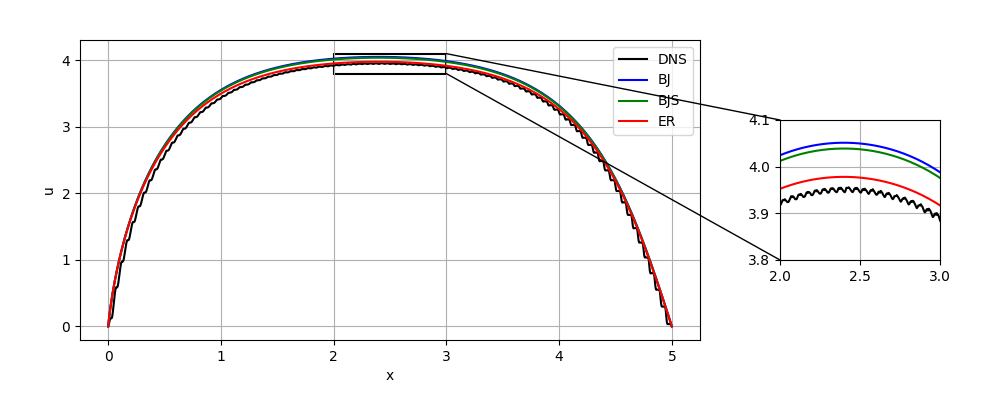}
    \includegraphics[width=0.4\linewidth]{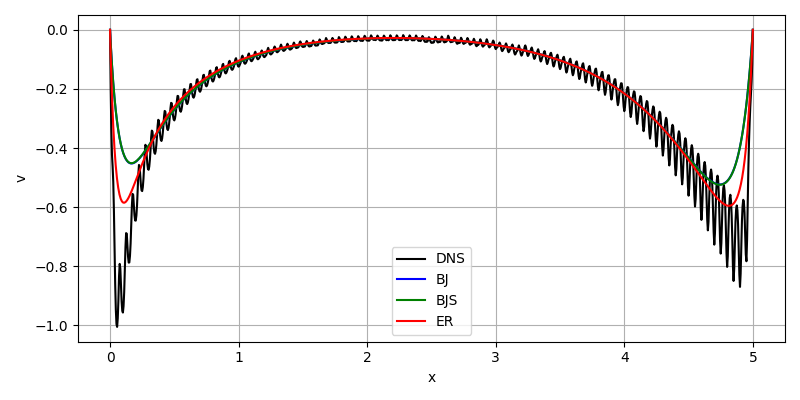}
    \caption{Tangential (left) and normal (right) velocity component along the fracture in \textit{Flow scenario~A} with geometry $G_5$ in $\Omega_1$ and geometry $G_2$ in $\Omega_2$.}
    \label{fig:Test2-G9G2}
\end{figure}

Now, to demonstrate that the new model is not only applicable to different porous media, but also to a vast range of flow scenarios, we consider two more test cases. As we already studied the influence of different pore structures on the flow profile, we present in the following only the results for the orthotropic geometry $G_4$ (Fig.~\ref{fig:permeability-table}).

\subsubsection{Flow scenario B}
\label{sec:B}
\begin{figure}[!ht]
    \centering
    \includegraphics[width=0.49\linewidth]{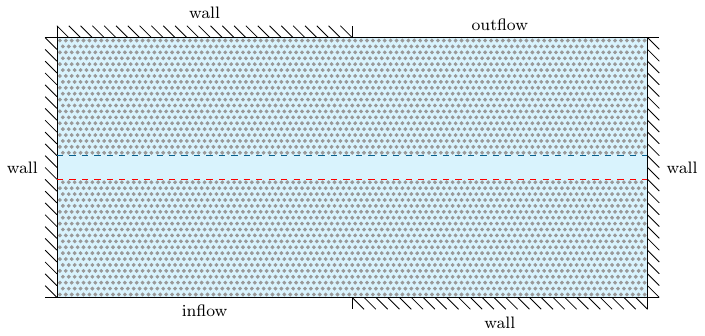}
    \raisebox{4.5mm}{
    \includegraphics[width=0.49\linewidth]{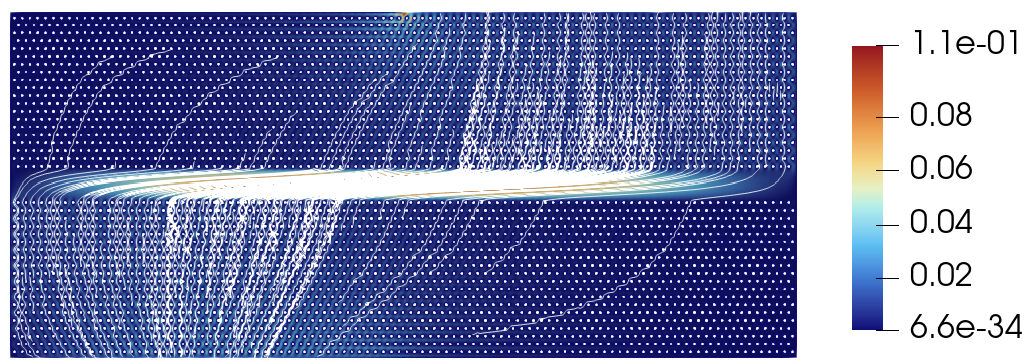}}
    \caption{Geometric setting and and pore-scale velocity magnitude in \textit{Flow scenario B.}
    }
    \label{fig:geometry-test3}
\end{figure}

In this flow scenario, we consider an inflow $\Gamma_2^{\mathrm{in}}$
 on the right bottom of the flow domain $\Omega$ and an outflow~$\Gamma_1^{\mathrm{out}}$ on the left top (Fig.~\ref{fig:geometry-test3}).
 We consider the flow domain $\Omega = [0, 5.0] \times [-1.1, 1.1]$, the inflow boundary $\Gamma_2^{\text{in}} = [0, 2.5] \times \{-1.1\}$, the outflow boundary $\Gamma_1^{\text{out}} = [2.5, 5.0] \times \{1.1\}$ and the remaining boundaries $\Gamma^{\text{wall}} = ([0, 2.5] \times \{1.1\}) \cup ([2.5, 5.0] \times \{-1.1\}) \cup (\{0\} \times [-1.1, 1.1]) \cup (\{5.0\} \times [-1.1, 1.1])$. In the hybrid-dimensional model, we have the porous-medium domains $\Omega_1 = [0,5.0] \times [0.1, 1.1]$ and $\Omega_2 = [0, 5.0] \times [-1.1, -0.1]$ and the complex interface $\Gamma = [0, 5.0] \times \{0\}$. We use the following boundary conditions
\begin{align}
    \vec{v} = \left(0,0.01\cdot x \cdot(x-2.5)\right) \text{ on } \Gamma_2^{\text{in}}, \quad p = 0 \text{ on } \Gamma_1^{\text{out}}, \quad \vec{v} = \textbf{0} \text{ on } \Gamma^{\text{wall}}.
\end{align}

\begin{figure}[!ht]
    \centering
    \includegraphics[width=0.55\linewidth]{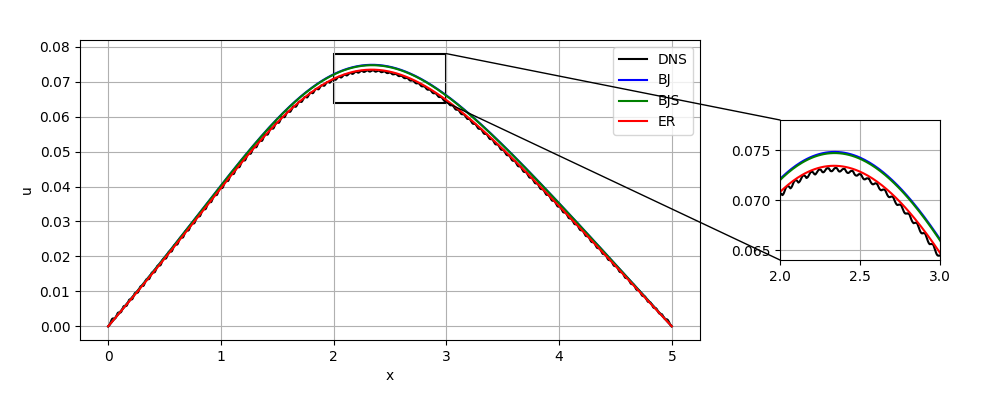}
    \includegraphics[width=0.4\linewidth]{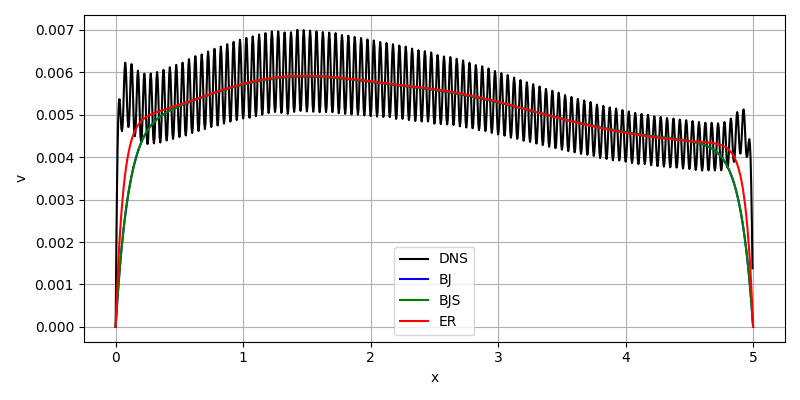}
    \caption{Tangential (left) and normal  (right) velocity component along the fracture in \textit{Flow scenario~B} with geometry $G_4$ in both porous-medium domains.}
    \label{fig:test3}
\end{figure}

As stated above, it is sufficient to consider one pore geometry, namely the orthotropic geometry $G_4$, for the remaining flow scenarios. Figure~\ref{fig:test3} compares the tangential and normal velocity profiles obtained with the hybrid-dimensional models to the corresponding pore-scale solution. Here, the newly developed hybrid-dimensional model also shows the best agreement with the pore-scale simulations in comparison to the models developed in~\cite{Rybak-Metzger-20}.

\subsubsection{Flow scenario C}
\label{sec:C}
In this flow scenario, we consider the flow domain $\Omega = [0, 5.0] \times [-1.1, 1.1]$, the inflow boundary $\Gamma_2^{\text{in}} = [0, 5.0] \times \{-1.1\}$, the outflow boundary $\Gamma_1^{\text{out}} = \{0\} \times [0.6, 1.1] \cup \{5.0\} \times [0.6, 1.1]$ and the remaining boundaries $\Gamma^{\text{wall}} = ([0, 5.0] \times \{1.1\}) ) \cup (\{0\} \times [-1.1, 0.6]) \cup (\{5.0\} \times [-1.1, 0.6])$. In the hybrid-dimensional model, we have the porous-medium domains $\Omega_1 = [0,5.0] \times [0.1, 1.1]$ and $\Omega_2 = [0, 5.0] \times [-1.1, -0.1]$ and the complex interface $\Gamma = [0, 5.0] \times \{0\}$. We use the following boundary conditions
\begin{align}
    \vec{v} = \left(0, 0.01\cdot x\cdot(x-5.0)\right) \text{ on } \Gamma_2^{\text{in}}, \quad p = 0 \text{ on } \Gamma_1^{\text{out}}, \quad \vec{v} = \textbf{0} \text{ on } \Gamma^{\text{wall}}.
\end{align}

\begin{figure}[!ht]
    \centering
    \includegraphics[width=0.49\linewidth]{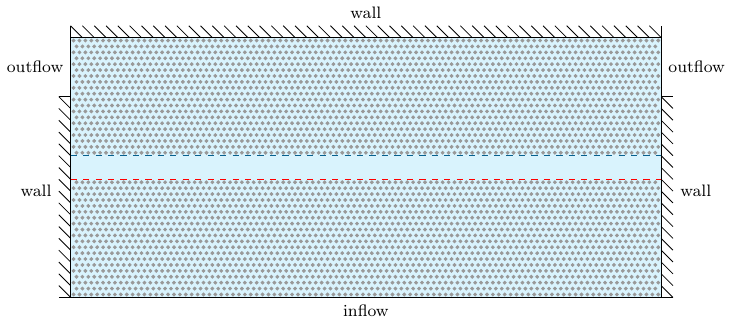}
    \raisebox{3.5mm}{
    \includegraphics[width=0.49\linewidth]{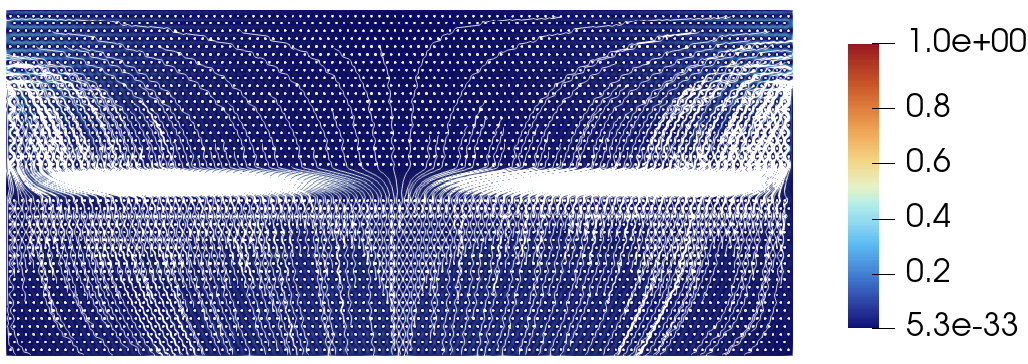}}
    \caption{Geometric setting and pore-scale velocity magnitude in \textit{Flow scenario C}.}
    \label{fig:geometry-test4}
\end{figure}

\begin{figure}[!ht]
    \centering
    \includegraphics[width=0.55\linewidth]{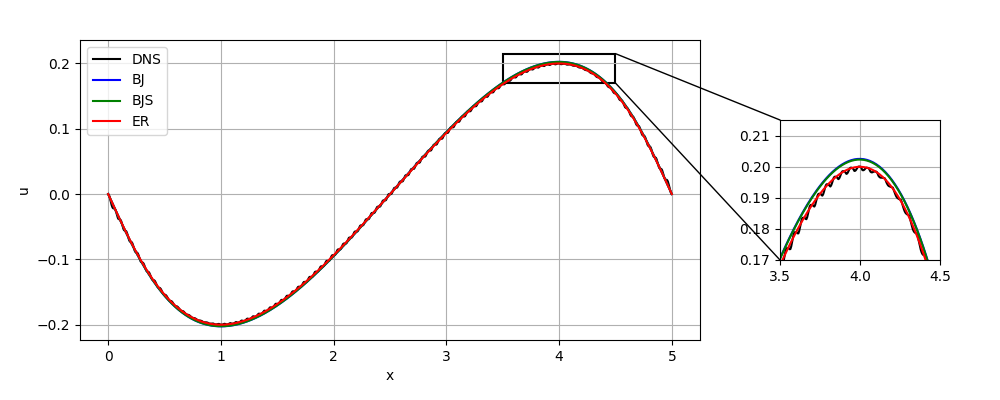} 
    \includegraphics[width=0.41\linewidth]{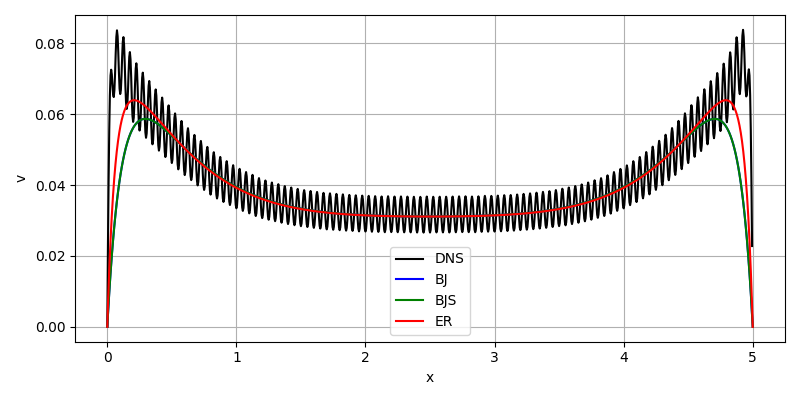}
    \caption{Tangential (left) and normal  (right) velocity component along the fracture in \textit{Flow scenario~C} with geometry $G_4$ in both porous-medium domains.}
    \label{fig:test4}
\end{figure}

For \textit{Flow scenario C}, similar observations can be made. As shown in Fig.~\ref{fig:test4}, the newly developed hybrid-dimensional model again demonstrates the best accuracy compared to the Stokes--Darcy hybrid-dimensional models available in the literature. These results confirm that the developed hybrid-dimensional model with generalised coupling conditions provides higher accuracy for a wide range of flow scenarios compared to the existing models.

\section{Conclusions}\label{sec:Conclusion}
In this work, we propose a new hybrid-dimensional Stokes--Darcy model to describe flow in fluid-filled fractured porous media. The hybrid-dimensional model is derived using vertical averaging from its full-dimensional counterpart. We consider the Stokes equations in the fracture surrounded by porous media with different material properties. The flow through the porous media is based on Darcy's law, and the generalised interface conditions are considered on the fracture-matrix interfaces. We prove the well-posedness of the derived hybrid-dimensional model for isotropic and orthotropic porous media. We validate the new hybrid-dimensional model using pore-scale simulations and compare it to previously developed hybrid-dimensional Stokes--Darcy models available in the literature considering a wide range of flow scenarios and pore geometries. The new model demonstrates the highest accuracy for all investigated test cases. We consider two directions for the future work. First, we will develop block preconditioners to accelerate numerical computations of the newly developed hybrid-dimensional model.  Second, we will compare our developed full-dimensional Stokes--Darcy fracture model to the TPM phase-field model from~\cite{rivas2025fluid, Rivas-etal-26} and validate our model with experimental data provided in these papers. 

\appendix
\renewcommand{\theequation}{A\arabic{equation}}
\setcounter{equation}{0}

\bibliographystyle{abbrv} 
\bibliography{references}

\end{document}